\documentclass{article}

\usepackage[hmargin={2.5cm,2.5cm}, vmargin={3cm,3cm}, dvips]{geometry}
\usepackage[utf8]{inputenc}
\usepackage{amssymb,amsmath,amsthm}
\usepackage{authblk}
\usepackage{graphicx}
\graphicspath{ {./img/} }
\usepackage{comment}
\usepackage[draft]{fixme}
\fxsetup{theme=color}
\usepackage{float}
\usepackage{enumitem}
\definecolor{fxtarget}{rgb}{0.8000,0.0000,0.0000}
\definecolor{fxnote}{rgb}{0.8000,0.0000,0.0000}
\usepackage{tikz}
\usetikzlibrary{calc,angles}
\usepackage{subcaption}
\usepackage[linesnumbered,ruled,vlined]{algorithm2e}

\usepackage{bm}
\usepackage{xcolor}

\newcommand\ls[1]{{\color{black}{#1}}}

\DeclareMathOperator*{\argmin}{arg\,min}

\newtheorem{theorem}{Theorem}[section]
\newtheorem{proposition}[theorem]{Proposition}
\newtheorem{corollary}[theorem]{Corollary}
\newtheorem{lemma}[theorem]{Lemma}
\theoremstyle{definition}
\newtheorem{definition}[theorem]{Definition}
\newtheorem{example}[theorem]{Example}

\theoremstyle{remark}
\newtheorem{remark}[theorem]{Remark}
\numberwithin{equation}{section}

\def \BC{{\mathcal C}}
\def \bis{{\bf b}}
\def \cross{{\bf n}}
\def\a{{\bf a}}
\def\b{{\bf b}}

\def\C{{\bf C}}

\def\S{{\bf S}}

\def\f{{\bf f}}

\def\h{{\bf h}}
\def\ii{{\bf i}}
\def\j{{\bf j}}
\def\k{{\bf k}}

\def\p{{\bf p}}
\def\q{{\bf q}}
\def\r{{\bf r}}
\def\s{{\bf s}}
\def\t{{\bf t}}
\def\u{{\bf u}}

\def\w{{\bf w}}
\def\vv{{\bf v}}
\def\X{{\bf X}}
\def\S{{\bf S}}
\def\A{{\mathcal A}}
\def\B{{\mathcal B}}

\def\Q{{\mathcal Q}}
 
\def\U{{\mathcal U}}
\def\W{{\mathcal W}}

\begin{document}

\title{Piecewise rational rotation-minimizing motions via data stream interpolation}

%    author one information

\author[$\dagger$]{Carlotta Giannelli}
\author[$\dagger$]{Lorenzo Sacco}
\author[$\dagger$]{Alessandra Sestini}
\author[$\star$]{Zbyn\v{e}k \v{S}\'{i}r}

\affil[$\dagger$]{Dipartimento di Matematica e Informatica ``U. Dini'',Universit\`{a} degli Studi di Firenze, Florence, Italy}
\affil[$\star$]{Charles University, Faculty of Mathematics and Physic, Prague, Czech Republic}

% \author{ Lorenzo Sacco, Alessandra Sestini, Zbyn\v{e}k \v{S}\'{i}r}
% \address{Dipartimento di Matematica e Informatica ``U. Dini'',Universit\`{a} degli Studi di Firenze, Florence, Italy}  
% \curraddr{}
% \email{carlotta.giannelli@unifi.it}
% \thanks{}

%    author two information
% \author{Lorenzo Sacco}
% \address{Dipartimento di Matematica e Informatica ``U. Dini'',Universit\`{a} degli Studi di Firenze, Florence, Italy}  
% \curraddr{}
% \email{lorenzo.sacco@unifi.it}
% \thanks{}

%    author three information
% \author{Alessandra Sestini}
% \address{Dipartimento di Matematica e Informatica ``U. Dini'',Universit\`{a} degli Studi di Firenze, Florence, Italy}  
% \curraddr{}
% \email{alessandra.sestini@unifi.it}
% \thanks{}

%    author four information
% \author{Zbyn\v{e}k \v{S}\'{i}r}
% \address{Charles University, Faculty of Mathematics and Physic, Prague, Czech Republic}  
% \curraddr{}
% \email{sir@karlin.mff.cuni.cz}
% \thanks{}

\date{}

% \dedicatory{}

\maketitle

%    Abstract is required.
\begin{abstract}
When a moving frame defined along a space curve is required to keep an axis aligned with the tangent direction of motion, the use of rotation-minimizing frames (RMF) avoids unnecessary rotations in the normal plane. The construction of rigid body motions using a specific subset of quintic curves with rational RMFs (RRMFs) is here considered. In particular, a novel geometric characterization of such subset enables the design of a local algorithm to  interpolate an assigned stream of positions, together with an initial frame orientation. To achieve this, the translational part of the motion is described by a parametric $G^1$ spline curve whose segments are quintic RRMFs,  with a globally continuous piecewise rational rotation-minimizing frame.  A selection of numerical experiments illustrates the performances of the proposed method on synthetic and arbitrary data streams.
%obtained from sampling but also generically  assigned to broadly design a reasonable space trajectory.
\end{abstract}

\begin{center}
\textbf{e-mail:} \\   carlotta.giannelli@unifi.it, lorenzo.sacco@unifi.it,\\  alessandra.sestini@unifi.it, sir@karlin.mff.cuni.cz
\end{center}


\section{Introduction}
In this paper we are interested in defining a special kind of rigid body motions based on piecewise rational parametric forms. The method follows a local approach and takes in input a stream of points, together with an assigned initial orientation of the rigid body. Any rigid body motion can be prescribed considering a parametric curve $ {\mathbf r}: [a\,,\,b] \rightarrow \mathbb{E}^3\,,\,\, \mathbf{r}{(t)}=(x(t),y(t),z(t)),$ that specifies the path described by the center of mass of the body, and an associated  triple $(\mathbf{f}_1(t),\mathbf{f}_2(t),\mathbf{f}_3(t))$  of mutually orthogonal unit vectors, that fixes the body orientation along $\mathbf{r}$.  The frame is usually denoted \emph{adapted} if,  at each parameter value, $\mathbf{f}_1(t)$  coincides with the unit tangent $\mathbf{t}(t)=\mathbf{r}'(t)/|\mathbf{r}'(t)|$ to the curve, while $\mathbf{f}_2(t),\mathbf{f}_3(t)$ span the normal plane orthogonal to $\mathbf{f}_1(t)$. To avoid unnecessary rotations in the normal plane, among the family of adapted frames, we consider motions where the frame is rotation-minimizing (RMF) \cite{bishop75}. As well as the Frenet frame, defined by $\t$, together with the normal \ls{$\mathbf{n}$} and binormal  \ls{$\mathbf{b}$} unit vectors of the curve $\mathbf{r}$, 
the RMF  is completely determined by $\mathbf{r}$, except for a constant rotation in the plane normal to the curve tangent. It was shown in \cite{guggenheimer89} that, when $(\mathbf{f}_1,\mathbf{f}_2,\mathbf{f}_3)$ is an RMF, then
$$
\mathbf{f}_2 = \cos\psi \ \mathbf{n} \,+\,  \sin\psi \ \mathbf{b}\,, \qquad
\mathbf{f}_3 = -\sin\psi \ \mathbf{n} \,+\,  \cos\psi \ \mathbf{b}\,,
$$
where $\psi = \psi(t)$ is the following angular function 
\begin{equation}\label{eq:guggh}
\psi(t) = \psi_0 - \int_0^t \tau({u}) \sigma(u) \mbox{d}u,
\end{equation}
with $\tau$ denoting the torsion of the curve and $\sigma = \vert \r'\vert$ its parametric speed. 
In general, the integral in \eqref{eq:guggh} does not admit an analytic \ls{representation among} standard polynomial and rational curves and related approximation schemes have been considered \cite{juettler99b,wang97,wang08}. 
As an alternative, in order to deal with (simple) exact analytic expressions, a suitable subset of (possibly piecewise) polynomial parametric curves can be considered so that the associated rotation-minimizing motions have \ls{a rational expression for $(\mathbf{f}_1,\mathbf{f}_2,\mathbf{f}_3)$}. For each spline segment, we then need to focus on a special subset of spatial polynomial Pythagorean--hodograph (PH) curves,  characterized by a rational unit tangent. In this paper we consider the subclass of PH curves with rational RMFs (RRMFs) introduced in \cite{fgms09}. We refer to \cite{review19} for a recent review of planar and spatial PH curves and related application algorithms. Note that any PH curve possesses a rational adapted frame, the so called Euler-Rodrigues frame (ERF), which admits a simple quaternion representation but is not necessarily a rotation--minimizing frame.

Several local algorithms based on polynomial PH curves have been proposed in the literature of the field for defining piecewise rational motions.  They usually require in input at least first order Hermite data and rely on the Euler-Rodrigues frame for prescribing the orientation component of the motion, see e.g., \cite{ch02,kv12,fhds13,ks21,kps24}. In particular, only the algorithm proposed in \cite{fhds13}, by relying on RRMF curves of degree seven, constructs an RMF motion. The algorithm proposed in this paper instead, relying on quintic RRMF curves whose RMF does not coincide with the ERF, produces the lowest degree RRMF spline curves. In addition, it requires in input just a  stream of points and a given initial frame orientation. A preliminary attempt to use RRMF quintic curves to interpolate first order Hermite data and prescribed end frame orientations was proposed in \cite{fgms2012}. It was there shown that it is not possible to guarantee the
existence of solutions in all cases, and a spline extension of the method was not considered. Note that the interpolation of end points and frames by rotation-minimizing motions strongly influences  the shape of the interpolants, see again \cite{fgms2012}.  

The focus of this paper is on the construction of a $G^1$ PH quintic spline curve, interpolating assigned positions at the spline knots and with a globally $C^0$ piecewise rational RMF. The algorithm is local, can be implemented in real-time, and  simultaneously produces for each spline segment the associated polynomial and rational forms of the curve and of the associated rotation minimizing frame, respectively. The input data for the construction of any spline segment are the initial position and first order Hermite information, the initial frame orientation, and the successive position to be interpolated. For solving the local Hermite interpolation problems here considered we rely on a special subset of RRMF quintics, usually denoted as RRMF quintics of class I \cite{fs12}. The development of our algorithm has been possible in view of a preliminary  analysis which presents a novel geometric characterization of this kind of curves. In particular, sufficient conditions for the existence of RRMF solutions to the considered local Hermite problem are detailed and properly exploited in the spline extension of the algorithm for an automatic generation of suitable unit tangents. %A unit tangent direction is automatically generated and interpolated in the spline extension of the algorithm. 

The remainder of the paper is organized as follows. After brieﬂy recalling the basic rules of the quaternion algebra, Section~\ref{sec:preliminaries} presents some properties and known results concerning PH and RRMF curves which are necessary for the following, along with some useful additional definitions. A new approach for the geometric characterization of RRMF quintic curves of class I, based on projections on the unit sphere, is introduced in Section~\ref{sec:characterization}, enabling the development of a new   local interpolation algorithm in Section~\ref{sec:local_alg}.  The $G^1$ continuous spline extension of the algorithm   is presented in Section~\ref{sec:global_alg}, together with some numerical examples. Finally, the key results of this work are summarized in Section~\ref{sec:conclusions}. 

\section{Preliminaries}\label{sec:preliminaries}
In this section some preliminary notions on spatial polynomial Pythagorean-hodograph (PH) curves are briefly summarized, with special focus on  the subset of PH curves whose rotation-minimizing frame is rational (RRMF curves). 
In order to make the paper self-contained, we start recalling the definition of polynomial PH curves.
\begin{definition}
A polynomial parametric curve $\mathbf r = \mathbf r(t)$ is a PH curve if its {\it parametric speed} $\sigma := \vert \mathbf r'\vert$ is a polynomial.
\end{definition}
Since we rely on the compact quaternion representation of spatial PH curves, the basic rules of quaternion algebra $\mathbb{H}$ are now reported. Any quaternion ${\A} \in \mathbb{H}$ can be defined as $(a_0,a_1,a_2,a_3)^T\,,$ with $a_i \in \mathbb R$. The notation $$a = scal(\A) := a_0 \qquad \text{and} \qquad \a = vect(\A) := (a_1,a_2,a_3)^T$$ refers to the scalar and vector part of the quaternion ${\A}, respectively.$ It is then possible to adopt the short representation ${\A} = a + \a$. If $a = 0\,, {\A}$ is called {\it pure vector} quaternion and can be shortly denoted just as $\a.$ The whole subset  of $\mathbb{H}$ of pure vector quaternions is denoted as $\mathbb{H}_v$ and identified with $\mathbb{R}^3$. Conversely, when $\a$ vanishes, ${\A}$ is a {\it pure scalar} quaternion and can just be denoted as any number in $\mathbb{R}$. The sum in $\mathbb{H}$ follows standard rules of sum in $\mathbb R^4$, quaternion product instead is non commutative and can be compactly defined as 
$$ {\A} {\mathcal B} = (a + \a)(b+\b) := (ab - \a \cdot  \b) + (a \b + b \a + \a \times \b)\,,$$
where conventional notation is used to denote scalar and cross vector products. ${\A}^* := a - \a$ denotes the conjugate of  ${\A},$ and, as a consequence, the product ${\A}{\A}^* = {\A}^*{\A} =a^2 + \a^T\a$ is a nonnegative pure scalar quaternion. The {\it module} $\vert {\A} \vert$ of a quaternion is defined as $\vert {\A} \vert := \sqrt{{\A}{\A}^*}$ and ${\A}$ is a {\it unit} quaternion if $\vert {\A} \vert = 1$. For any pure vector $\vv$, it is also true  that the  product  ${\A}\,\vv\, {\A}^*$  remains a pure vector quaternion. Keeping this in mind, we can easily understand how unit quaternions allow a compact representation of spatial rotations. Indeed, since any unit quaternion $\U$ can be  represented in terms of an angle $\varphi$ and of a unit vector $\ii$ as $\U = e^{\ii\varphi} :=  \cos(\varphi) + \ii \sin(\varphi),$  with some quaternion computations it can be verified that, for any pure vector quaternion $\vv$, it is always true that  the vector 

\begin{equation} \label{rotation}
      \w = \U \, \vv \, \U^*, \qquad \text{with} \qquad \U = e^{\ii \left(\theta/2 \right)} 
\end{equation}
is obtained by rotating $\vv$ through the angle $\theta$ about the axis defined by $\ii$. The quaternion algebra offers  a compact representation of spatial PH curves, because of the following result \cite{choi02,farouki02}.

\begin{proposition}
   A polynomial parametric curve $\mathbf r$ is a PH \ls{curve} if and only if there exist a quaternion pre-image polynomial $\A\in \mathbb{H}[t]$ and a real polynomial $\rho\in \mathbb{R}[t]$ with no odd-multiplicity real root such that its hodograph $\mathbf h := \mathbf r'$ has the following form,
   \begin{equation}\label{PHrep}
       \mathbf h = \rho\, \A\, \ii\, \A^*,     
   \end{equation}
where $\ii$ denotes any unit pure vector.     
\end{proposition}
Without loss of generality, we can assume $\rho$ positive. Note that it is necessary to require that $\rho$ has no odd-multiplicity real root to guarantee a polynomial analytical expression for $\sigma = \vert \r' \vert = \vert \rho \vert \ \A\A^*$ in the whole $\mathbb{R},$  even if in practice it could be sufficient to require this assumption in the   parameter domain of the curve, let us say without loss of generality $[0,1].$ Since it will be useful in the following, we also define a special subset of PH curves.
\begin{definition}
A polynomial PH curve $\mathbf r$   is degenerate if its pre-image $\A$ has at least one real root belonging to its  parameter domain $[0,1].$
\end{definition}
If a PH curve is degenerate, there exists   $\Q \in \mathbb{H}[t]$ never vanishing in $[0,1]$ and $\chi \in \mathbb{R}[t]$ such that 
\begin{equation} \label{hoddegener}
\mathbf h = \rho\, \chi^2 \Q\, \ii\, \Q^* \,.
\end{equation}
Among many advantages, PH curves allow the exact computation of geometric quantities such as the arc length. Several interpolation algorithms have been developed in the latest years, see e.g.,  \ls{\cite{sj2007,gss22,kps24}}. Another interesting feature of PH curves is their rational tangent indicatrix. We recall that for any regular parametric curve $\r$ (possibly also with loops), the tangent indicatrix is the locus of points on the unit sphere $\mathbb{S}^2$ traced by its normalized tangent vectors, interpreted as spherical points.
From now on, we will equivalently refer to spherical points and to  corresponding unit vectors in $\mathbb{R}^3$ with origin in the center of the sphere.

\
\begin{corollary}
    The tangent indicatrix $\mathbf t$ of a non degenerate polynomial PH curve with hodograph as in \eqref{PHrep} has the rational parametric form 
    \begin{equation}\label{rtind}
        \mathbf t = \mathbf t(t) = \frac{\mathcal A(t)\, \mathbf i \,\mathcal A(t)^*}{\mathcal A(t)\mathcal A(t)^*}\,, \qquad t \in [0,1].
    \end{equation}
    
\end{corollary}
Note that the tangent indicatrix of a degenerate PH curve has the same compact expression with $\Q$ replacing $\A$ and, in any case, the rational form of $\mathbf t$  does not depend on $\rho.$ 
\begin{remark} \rm \label{rmk:ro}
Considering in \eqref{PHrep} a non zero degree real polynomial $\rho$ with all odd-multiplicity real roots not belonging to $[0,1]$ is safe from the point of view of cusps and it could be useful to gain some more flexibility with low complexity increase. This possibility, however, is not taken into account in this paper.    
\end{remark}

This work often presents pure vector quaternion equations whose solution is characterized by one angular degree of freedom as, for example, it is stated in the following remark.
\begin{remark}\label{fibres}\rm
For any choice of the unit vector quaternion $\ii,$ 
$e^{\ii\,\alpha}\,
\ii\,e^{-\ii\,\alpha} = \ii,$ as for standard complex numbers.
Consequently, the hodograph representation in \eqref{PHrep} does not change if the quaternion polynomial $\A = \A(t)$ is replaced by $\tilde \A(t)=\A(t)\,e^{\ii\,\alpha}.$  This implies that any construction of PH curves leads to a one-parameter system of pre-image solutions, all representing the same hodograph. 
\end{remark}

The characterization of RRMF curves was thoroughly studied over the years. In quaternion form it can be expressed as follows, see e.g., \cite{han08} and also \cite{fggss17} for additional theoretical insights.
\begin{proposition}    
A non degenerate polynomial PH curve with pre-image $\A$ has a rational RMF if and only if there exists a quaternion polynomial $\W = a + \ii\, b\,, \, a,b \in \mathbb{R}[t]$ s.t.
\begin{equation} 
        \frac{\it{scal}(\A' \, \ii \,\A^*)}{\A \, \A^*} = \frac{\text{scal}(\W' \, \ii\, \W^*)}{\W \, \W^*}\,.
        \label{han08}
\end{equation}
\end{proposition}

We observe that the easiest way to fulfill \eqref{han08} consists in selecting $b = 0$, which implies $\text{scal}(\W' \, \ii\, \W^*) = 0$, and simultaneously requiring  $\text{scal}(\A' \, \ii \A^*) \equiv 0$. This corresponds to impose that the RMF $(\mathbf{f}_1,\mathbf{f}_2,\mathbf{f}_3)$ coincides with the so called Euler-Rodriguez frame (ERF) $(\mathbf{e}_1,\mathbf{e}_2,\mathbf{e}_3)$, where
\begin{equation*}
    \mathbf{e}_1 := \frac{\A\,\ii\,\A^*}{\A \A^*},\qquad 
    \mathbf{e}_2 := \frac{\A\,\j\,\A^*}{\A \A^*}, \qquad 
    \mathbf{e}_3 := \frac{\A\,\k\,\A^*}{\A \A^*}\,.
\end{equation*}
which is clearly rational. By defining 
$$\B := \A\, \W\,,$$
the RMF can be expressed in the following compact form \cite{fgms2012} 
\begin{equation*}
    \mathbf{f}_1 := \frac{\B\,\ii\,\B^*}{\B\,\B^*} = \frac{\A\,\ii\,\A^*}{\A\, \A^*},\qquad 
    \mathbf{f}_2 := \frac{\B\,\j\,\B^*}{\B\, \B^*}, \qquad 
    \mathbf{f}_3 := \frac{\B\,\k\,\B^*}{\B\, \B^*}.
\end{equation*}
Note that, in general, the rational form of the frame is of higher degree than the one of the tangent indicatrix of the curve.

\subsection{PH quintics}
Let us now review the characterization of quintic PH curves, the class of PH curves considered in our algorithms. In particular, we will study the family of PH quintics which share the same tangent indicatrix,  
\ls{considering B\'ezier representations, defined in terms of Bernstein polynomials 
\[
B_i^n(t) := \binom{n}{i}\,t^i\,(1-t)^{n-i}\,, \quad i = 0,\ldots,n\,.
\]}

Given a second degree pre-image in B\'ezier form, 
    \begin{equation}\label{preim}
        \mathcal A(t)=\mathcal A_0 B^2_0(t)+\mathcal A_1 B^2_1(t)+\mathcal A_2 B^2_2(t)\,, \quad t \in [0,1]\,,
    \end{equation}
      with quaternion coefficients $\mathcal A_i, i=0,1,2$, the associated hodograph is the following parametric curve, 
    \begin{equation} \label{preimage}
    \mathbf h(t)=\mathcal A(t)\,\mathbf i\, \mathcal A^*(t)=
    \sum_{i=0}^4\mathbf h_i\, B^4_i(t)\,, \quad t \in [0,1]\,,
    \end{equation}
   where 
\begin{align}\label{hod}
    &\mathbf h_0:= \mathcal A_0\,\mathbf i\,
    \mathcal A^*_0, 
    \qquad 
    \mathbf h_1:= \frac{1}{2}
    (\mathcal A_1\, \mathbf i\,\mathcal A^*_0+\mathcal A_0\,
    \ii\, \mathcal A^*_1),\qquad 
    \h_2:= \frac{1}{6}
    (\mathcal A_2\, \mathbf i\,\mathcal A^*_0+4\, \mathcal A_1\, \mathbf i\,\mathcal A^*_1+\mathcal A_0\,\mathbf i\,
    \mathcal A^*_2),\nonumber\\
    &
    \h_3:= \frac{1}{2}
    (\A_2\, \mathbf i \,\mathcal A^*_1+\mathcal A_1\,
    \mathbf i \,\mathcal A^*_2),
    \qquad \mathbf h_4 := \mathcal A_2\,
    \mathbf i\,\mathcal A^*_2\,.
\end{align}
Integrating \eqref{preimage} gives the B\'ezier form
\begin{equation}\label{PHcurve}
    \r(t) = \sum_{i=0}^5\mathbf r_i\, B^5_i(t)\,, \quad t \in [0,1]\,,
\end{equation}
of the PH quintic, with $\mathbf r_0$ denoting an integration constant in \ls{the affine space} $\mathbb{E}^3$ and the other control points $\r_i, i=1,\ldots,4,$ defined in terms of $\mathcal A_i, i=0,1,2$ and of  $\mathbf r_0$ as follows,
\begin{align}\label{PHcontrolpoints}
    &\mathbf r_1:= \r_0 + \frac{1}{5}\, \mathcal A_0\,\mathbf i\,
    \mathcal A^*_0, 
    \,\,
    \mathbf r_2 := \r_1 + \frac{1}{10}
    (\mathcal A_1 \mathbf i\mathcal A^*_0+\mathcal A_0
    \ii \mathcal A^*_1),
    \,\,
    \r_3 := \r_2 + \frac{1}{30}
    (\mathcal A_2\, \mathbf i\,\mathcal A^*_0 +4\, \mathcal A_1\, \mathbf i\,\mathcal A^*_1+\mathcal A_0\,\mathbf i\, 
    \mathcal A^*_2), \nonumber\\
    &\r_4 := \r_3 + \frac{1}{10}
    (\A_2 \mathbf i\mathcal A^*_1+\mathcal A_1
    \mathbf i\mathcal A^*_2),
    \,\,
    \mathbf r_5 := \r_4 + \frac{1}{5}\,\mathcal A_2\,
    \mathbf i\,\mathcal A^*_2.
\end{align}

Considering for simplicity only non degenerate PH quintics, the associated tangent indicatrix $\mathbf t$ is a curve on the unit sphere with the following degree $4$ rational parametric form,
\begin{equation} \label{ind}
    \mathbf t(t) = \frac{\mathcal A(t)\,\mathbf i\, \mathcal A^*(t)}{\mathcal A(t)\,\mathcal A^*(t)}
   =\frac{\sum_{i=0}^4\h_i\,B^4_i(t)}{\sum_{i=0}^4 w_i\,B^4_i(t)}
\end{equation}
with  weights 
\begin{align}\label{weights}
    &w_0:= \mathcal A_0
   \mathcal A^*_0\,,
   \qquad
   w_1:= \frac{1}{2}\,
   (\mathcal A_1\, \mathcal A^*_0+\mathcal A_0\,
   \mathcal A^*_1)\,,
   \qquad
   w_2:= \frac{1}{6}\,
   (\mathcal A_2\, \mathcal A^*_0+4 \,\mathcal A_1 \,\mathcal A^*_1+\mathcal A_0\,
   \mathcal A^*_2)\,,\nonumber\\
   &w_3:= \frac{1}{2}\,
   (\mathcal A_2\, \mathcal A^*_1+\mathcal A_1\,
   \mathcal A^*_2)\,,
   \qquad
    w_4:= \mathcal A_2\,
   \mathcal A^*_2\,.
\end{align}
In principle, the weights $w_i$, for $i=1,2,3$, can be negative, but the denominator of \eqref{ind} is always a positive polynomial. When the weights are all not vanishing, we may express the tangent indicatrix in the following rational B\'ezier form
\begin{equation} \label{ind2}
    \mathbf t(t)=\frac{\sum_{i=0}^4 w_i \,\t_i\, B^4_i(t)}{\sum_{i=0}^4 w_i\, B^4_i(t)}
\end{equation}
where 
\begin{equation}\label{pts}
    \t_i=\mathbf h_i/w_i\,, \quad w_i \ne 0.
\end{equation}

We now recall a general result about the weights of rational B\'ezier curves, see, e.g.,  \cite{hoschek93}. 

\begin{lemma}\label{rep}
Let 
\[
{\mathbf r}(t)=
\frac{\sum_{i=0}^n w_i\, \mathbf r_i \, B^n_i(t)}{\sum_{i=0}^n\, w_i B^n_i(t)}
\] 
be a rational Bézier curve of degree $n$ with control points $(\r_0, \ldots, \mathbf \r_n)$ and weights $(w_0, \ldots, w_n)$. Let  $\lambda \in \mathbb R^+$ be a positive real number and define $\tilde {\mathbf r}(\tilde t\,):=\mathbf r(\xi(\tilde t \,))$ as the reparameterization of ${\mathbf r}(t)$ via the linear rational reparameterization
    \begin{equation}\label{repar}
    \xi(\tilde t)=\frac{\lambda \tilde t}{(\lambda-1)\tilde t+1}.
\end{equation} 
Then $\tilde {\mathbf r}(\tilde t \,)$ is a rational B\'ezier curve with the same control points $(\mathbf r_0, \ldots, \mathbf r_n)$ and new weights $(\tilde w_0, \ldots, \tilde w_n)$, where $\tilde w_i = \lambda^i w_i$. 
\end{lemma}

The following proposition about the tangent indicatrix $\t$ can then be proved. 

\begin{proposition}
    \label{PropRepar}
    Let us have two PH curves with the following two quadratic pre-images 
    \begin{equation*}
    \mathcal A(t) = \mathcal A_0 B^2_0(t)+\mathcal A_1 B^2_1(t)+\mathcal A_2 B^2_2(t), \qquad
    \tilde {\mathcal A}(\tilde t)=\tilde {\mathcal A}_0 B^2_0(\tilde t)+\tilde {\mathcal A}_1 B^2_1(\tilde t)+\tilde {\mathcal A}_2 B^2_2(\tilde t)        
    \end{equation*}
    where $$\tilde{\mathcal A}_0=\mu\, \mathcal A_0, \quad \tilde{\mathcal A}_1=\mu\,\lambda \mathcal A_1, \quad\tilde{\mathcal A}_2=\mu \,\lambda^2 \mathcal A_2$$
    for some positive real numbers $\mu$, $\lambda$. Then for their tangent indicatrices $\mathbf t$ and $\tilde{\mathbf t}$ it holds 
    $\tilde{\mathbf t}(\tilde t\,)=\mathbf t(\xi(\tilde t\,)),$ with $\xi: [0,1] \rightarrow [0,1]$ defined \ls{by} \eqref{repar}.
Consequently, $\t$ and $\tilde \t$ have the same image on the sphere.
\end{proposition}
\begin{proof}
    From the formulas \eqref{hod} and \eqref{weights} we see that $\tilde {\mathbf h}_i=\mu^2\lambda^i\, \mathbf h_i$ and $\tilde w_i=\mu^2\lambda^i w_i$. Equation \eqref{pts} thus implies $\tilde \t_i= \t_i$ and the two curves thus have the same control points. Concerning the weights the common factor $\mu^2$ can be disregarded because it cancels out in \eqref{ind} and the weights are as in Lemma \ref{rep} in view of the considered reparameterization.
\end{proof}

\begin{figure}[t!]
\begin{center}  \includegraphics[width=0.33\textwidth]{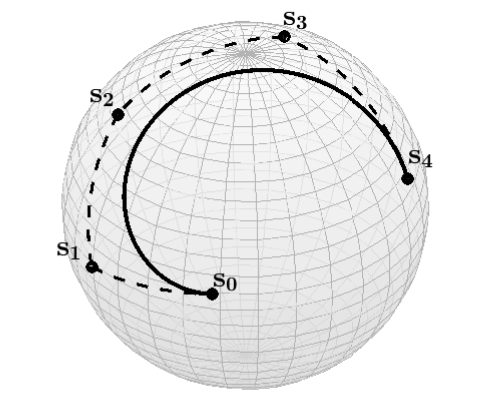}
    \caption{Example of a degree 4 tangent indicatrix (solid black line) of a regular PH quintic with its spherical control polygon (dashed black line) and  spherical control points $\s_i, i=0,\ldots,4$ (black dots).}
    \label{fig:spherical}
\end{center}
\end{figure}

If $\h_i \ne {\bf 0}$ for all $i=0,\ldots,4,$ we can associate to the PH quintic introduced in \eqref{PHcurve} the set of {\em spherical control points} belonging to $\mathbb{S}^2$ 
\begin{equation}\label{sphericalCP}
    \s_i := \frac{\mathbf h_i}{\vert\mathbf h_i\vert}, \quad i=0,\ldots,4.
\end{equation}
Without being real control points of $\mathbf t$, as shown for  example in Figure~\ref{fig:spherical}, these spherical points will be very useful both for explaining the geometry of the hodograph in Section~\ref{sec:characterization} and for explaining our local Hermite interpolation approach in Section~\ref{sec:local_alg}.

Concerning PH curves of degree less or equal to five, we remind that non planar RRMFs with RMF coincident with the ERF do not exist, see \cite{ch02}. Furthermore, relying on the Hopf map representation of PH curves, the case of RRMF cubics and quintics whose RMF is not the ERF has been deeply studied in \cite{fgms09}, where the PH subset fulfilling condition \eqref{han08} with non real $\W$ and deg$(\W ) =$ deg$(\A) =2$ was determined. Subsequently, the following compact characterization in quaternion form was obtained \cite{f2010}. 

\begin{proposition}  \label{charRRMF}
A PH quintic with the quadratic quaternion preimage in B\'ezier form \eqref{preim} verifies \eqref{han08} with deg$(\W) = $ deg$(\A) = 2$ and vect$(\W) \neq \bm{0}$, if and only if   
\begin{equation} \label{rrmf}
    \A_1\,\ii\, \A_1^* = vect(\A_2 \,\ii\, \A_0^*).
\end{equation}
\end{proposition}
Note that planar PH quintics (which are surely RRMF curves) do not necessarily verify the condition in \eqref{rrmf} but they still fulfill condition \eqref{han08} with $\W$ real. We conclude this subsection observing that in \cite{fs12} the set of quintic RRMFs characterized by the requirement in (\ref{rrmf}) was called class I quintic RRMFs, while those whose RMF coincides with the ERF were said of class III. The other type of quintic RRMFs there studied were those of class II corresponding to assuming deg$(\W) =1.$ These are not necessarily planar curves and so they could be here of interest as well. However, as already mentioned in \cite{fs12}, see also \cite{farouki16}, the available description of such subset seems less suited for deriving geometric constructive algorithms. For this reason in this paper  we rely just on quintic RRMFs of class I. 

%---------------------
\subsection{Vector and quaternion formulas}

We now present some definitions that will be useful to simplify the presentation of our results.

\begin{definition} 
Given two quaternions $\A$, $\B,$ let
\begin{equation} \label{com_mult} 
    \A \,\star\, \B := \frac{1}{2}\left( \A\, \ii\, \B^* + \B\, \ii\, \A^* \right), \qquad
\A \,\square\, \B := \frac{1}{2}(\A\, \B^* - \B\,\A^*),
\end{equation}
be two binary operators acting from $\mathbb{H} \times \mathbb{H}$ to $\mathbb{H}_v$.
\end{definition}

\begin{definition}
Given two spherical point $\mathbf v$, $\mathbf w,$   their minimum distance on $\mathbb{S}^2$ (geodetic distance) is denoted as angular distance and it can be identified as the smaller angle $\gamma \in [0,\pi]$ between them, so that $\cos(\gamma) = \mathbf v \cdot \w.$ 
\end{definition}

\begin{definition}
Given two vectors $\mathbf v$, $\mathbf w$,   their unit bisector is defined as    
\begin{equation} \label{bdef}
\bis(\mathbf v,\mathbf w):=\frac{\frac{\mathbf v}{\vert {\mathbf v} \vert}+ \frac{\mathbf w}{\vert{\mathbf w}\vert}}{\left \vert\frac{\mathbf v}{\vert{\mathbf v}\vert}+ \frac{\mathbf w}{\vert{\mathbf w}\vert}\right \vert}, \qquad 
\text{if}\qquad 
\frac{\mathbf v}{\vert \mathbf v \vert} +\frac{\w}{\vert \w \vert} \neq \mathbf{0}, 
\end{equation} 
and their negatively oriented normalized cross product as     
\begin{equation} \label{ndef}
\cross(\mathbf v, \mathbf w):= -\frac{\mathbf v \times \w}{\left \vert \mathbf v \times \w \right \vert}, \qquad \text{if}\qquad \mathbf v \times \w \neq \mathbf{0}.
\end{equation}
Note that $\mathbf b$ and $\mathbf n$ can be interpreted also as spherical points on $\mathbb{S}^2$.
\end{definition}

\begin{definition}\label{BC}
For any two spherical points $\mathbf v$ and $\w,$ the great circle of points on $\mathbb{S}^2$ having the same angular distance from  $\mathbf v$ and $\w$ is denoted as $\BC(\mathbf v,\w)$. If $\mathbf v \times \mathbf w \ne \mathbf{0},$ the order of the arguments $\mathbf v$ and $\w$ in $\BC(\mathbf v,\w)$ is used to specify that its orientation is given by the direction of the unit vector $\cross(\mathbf v,\w)$, tangent to $\BC(\mathbf v,\w)$ at the spherical point $\bis(\mathbf v,\w)$.
\end{definition}

\begin{figure}[t!]
\begin{center}
    \includegraphics[width=0.33\textwidth]{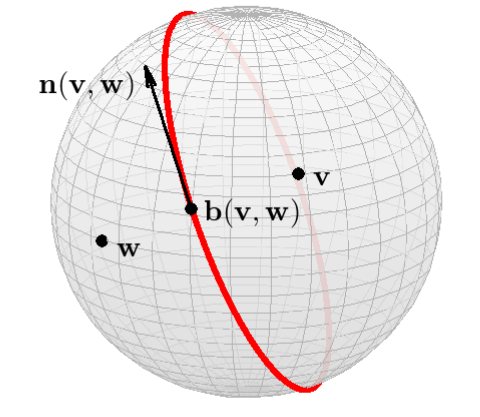}
    \caption{Two unit vectors $\vv$, $\w$ and their unit bisector $\bis(\vv,\w)$ all depicted as spherical points (black dots) on the unit sphere, together with the great circle $\BC(\vv,\w)$ (red line) and the negatively oriented normalized cross product $\cross(\vv, \w)$, represented as a vector (black arrow) applied at the point $\bis(\vv,\w)$.}
    \label{fig:bc}
\end{center}
\end{figure}

See also Figure~\ref{fig:bc} for the notation introduced in the previous two definitions. 

\begin{definition}
For any unit pure vector $\ii$ and any non zero pure vector $\mathbf{v}$, the quaternion square roots of $\mathbf{v}$ with respect to $\ii$ are all the quaternions $\A$ which are solutions of the quadratic vector equation 
\[
    \A\,\ii\,\A^* = \mathbf{v}\,.
\]
\end{definition}
\noindent
In any case the set of the quaternion roots of $\mathbf v$ with respect to $\ii$ is a one-parameter family of quaternions defined in the general setting as 
\begin{equation} \label{sqrt_quaternion}
    \A = \sqrt{\vert \mathbf v \vert} \,\, \bis(\ii,\mathbf{v})\ e^{\ii \alpha}\,, \qquad \mbox{if } \mathbf v \ne - \vert \mathbf v \vert \ \ii,
\end{equation}
and just as $A = \sqrt{\vert \mathbf v \vert}(\bm{\delta}^\perp_1 \cos \alpha+ \bm{\delta}^\perp_2 \sin \alpha)\,,$ in the special case $\mathbf v = - \vert \mathbf v \vert\, \ii\,,$
where $(\mathbf v/\vert \mathbf v \vert\,,\,\bm{\delta}^\perp_1\,,\,\bm{\delta}^\perp_2)$ is an orthonormal base \ls{of} $\mathbb R^3$.
We can observe that in equation~\eqref{sqrt_quaternion} the module of $\A$ is immediately identified as $\sqrt{\vert\vv\vert}$ and that $\U = \A / \vert \A \vert$ can be interpreted as a unit quaternion defining a rotation of $\ii$ into $\vv / \vert \vv \vert$, see \eqref{rotation}. In particular, if $\alpha = 0$, the rotation is computed around the bisector $\bis(\ii,\vv)$ through an angle $\theta = \pi$.

\section{Analysis of RRMF quintics of class I}
\label{sec:characterization}
In this section we analyze the hodographs of RRMF quintic curves of class I (RRMF5-I), whose pre-image is characterized by the algebraic condition \eqref{rrmf}. In particular, sufficient and necessary conditions are studied for both the hodograph control points $\mathbf h_i$ and their spherical versions $\s_i=\mathbf h_i/|\mathbf h_i|, i=0,\ldots,4,$ to generate RRMF5-I curves. 
Using the binary star operator introduced in \eqref{com_mult}, it is useful to rewrite the characteristic equation (\ref{rrmf}) as
\begin{equation}\label{rrmfstar}
   \A_1 \, \star \, \A_1 = \A_0 \, \star \, \A_2\,.
\end{equation}
This also allows us to rewrite the hodograph control points introduced in \eqref{hod}  as
\begin{equation}\label{control_point_rrmf}
\begin{array}{rcl}
    &\h_0= \A_0 \, \star \, \A_0,
   \quad
   \h_1= \A_0 \, \star \, \A_1,
   \quad
   \h_2= \A_1 \, \star \, \A_1 = \A_0 \, \star \, \A_2,
   \quad
   \h_3= \A_1 \, \star \, \A_2,
   \quad
   \h_4= \A_2 \, \star \, \A_2.
\end{array}
\end{equation}
For clarity of presentation,  we introduce the principal results together with some examples in Section~\ref{sec:gc}, and then present the related proofs in Subsection~\ref{sec::proof}.

\subsection{Geometric characterization}\label{sec:gc}
We start with a simple necessary condition for the spherical control points.

\begin{proposition}\label{pr1}
        Let $\s_i$, $i=0,\ldots,4$, be the spherical control points of the hodograph of an RRMF5-I curve. Then $$\mathbf s_1\in \BC(\mathbf s_0,\mathbf s_2), \qquad\mathbf s_2\in \BC(\mathbf s_0,\mathbf s_4),\qquad  \mathbf s_3\in \BC(\mathbf s_4,\mathbf s_2),
      $$ where \(\BC\) is introduced in Definition~\ref{BC}.
\end{proposition}
\begin{figure}[t!]
\begin{center}

    \includegraphics[width=0.329\linewidth]{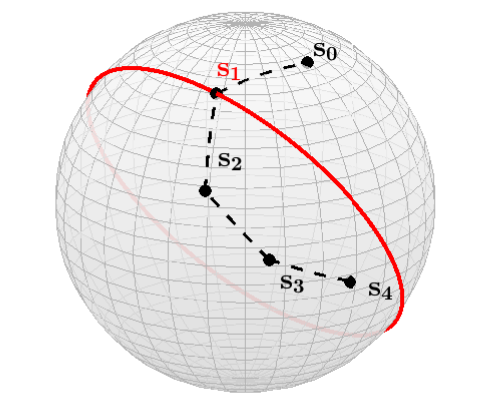}
    \includegraphics[width=0.329\linewidth]{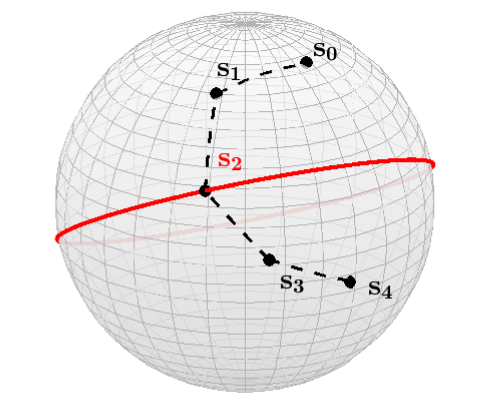}
    \includegraphics[width=0.329\linewidth]{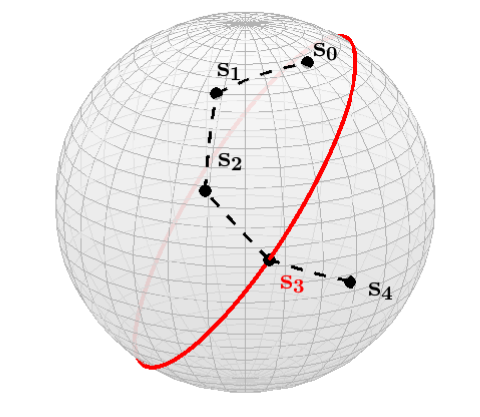}
    \caption{The spherical points $\s_0, \ldots, \s_4$ (black dots) of an RRMF5-I curve and the corresponding spherical control polygon (dashed black line). \ls{Each great circle (red lines) where $\s_1$ (left), $\s_2$ (center), and $\s_3$ (right) \ls{are located is} also shown.}}
    \label{fig:prop1_control_points}    
\end{center}
\end{figure}

A graphical representation of Proposition~\ref{pr1} is presented in Figure~\ref{fig:prop1_control_points}. Note that the order of the arguments in $\BC(\s_4,\s_2)$ is  motivated by the orientation of the great circles, as clarified by the following proposition, which also provides a sufficient condition for the spherical control points of an RRMF5-I curve. 

\begin{proposition}\label{pr2}
    Let three spherical points $\s_0$, $\s_4$, and $\s_2\in \BC(\s_0,\s_4)$ be given. Then there is a one-parameter family of pairs of spherical points $\s_1\in \BC(\s_0,\s_2)$ and $\s_3\in \BC(\s_4,\s_2)$ so that $\s_0$, $\s_1$, $\s_2$, $\s_3$, $\s_4$ are spherical control points of an RRMF5-I curve. Moreover, any point of $\BC(\s_0,\s_2)$ and $\BC(\s_4,\s_2)$ occurs precisely once in this way and the bijection $\s_1 \leftrightarrow \s_3$ between the two great circles is smooth and preserves the orientation. 
\end{proposition}
\begin{figure}[t!]
    \includegraphics[width=0.329\textwidth]{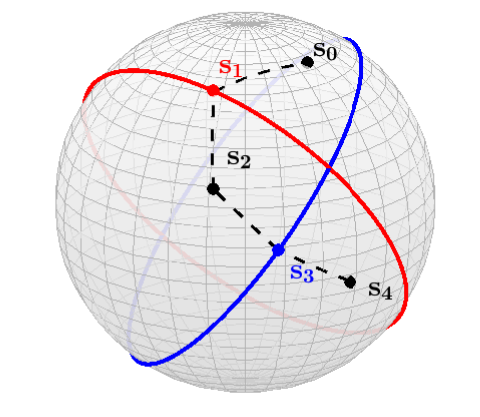}
    \includegraphics[width=0.329\textwidth]{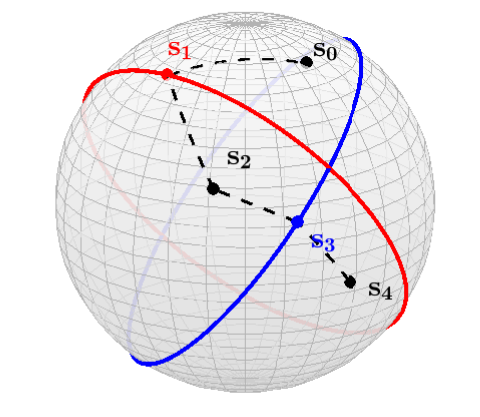}
    \includegraphics[width=0.329\textwidth]{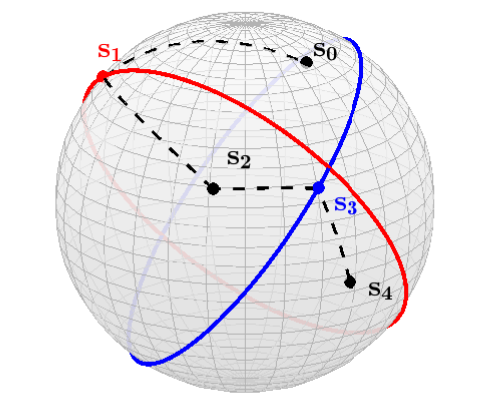}
    \caption{Three fixed spherical points $\s_0$, $\s_4$ and $\s_2\in\BC(\s_0,\s_4)$ (black dots), together with three possible positions (left, center, right) of the spherical point $\s_1$ \text{(red dots)} $\in \BC(\s_0,\s_2)$ \text{(red lines)} and corresponding positions of $\s_3$ \text{(blue dots)} $\in \BC(\s_0,\s_2)$ \text{(blue lines)} characterizing an RRMF5-I. The corresponding spherical control polygon (dashed black line) is also shown. }
    \label{fig:prop2_control_points}
\end{figure}

A graphical representation of Proposition~\ref{pr1} is presented in Figure~\ref{fig:prop2_control_points}. The explicit computation of the bijection $\BC(\s_0,\s_2)\leftrightarrow \BC(\s_4,\s_2)$ is quite technical and it will be provided in Remark~\ref{bijection}. Since $\mathbf h_i=\vert\mathbf h_i\vert\s_i$, for $i=0,\ldots,4$, to fully describe the geometry of all RRMF quintics, we need to specify not only the position $\s_i$ but also the lengths $\vert\mathbf h_i\vert$.

\begin{proposition}
\label{pr3}
Let $\s_0, \ldots, \s_4$ be an admissible configuration of spherical points as described in Proposition \ref{pr2}. For any choice of nonzero lengths $\vert\mathbf h_0\vert$ and $\vert \mathbf h_4\vert$, there exists precisely one set of control points $\mathbf h_0, \ldots, \mathbf h_4$ representing the hodograph of an RRMF5-I. Moreover, for various choices of $\vert\mathbf h_0\vert$ and $\vert \mathbf h_4\vert$, the corresponding tangent indicateces are related via a linear rational reparameterization. Consequently, the image of the tangent indicatrix is fully characterized by the points $\s_i$, $i=0,\ldots,4$.
\end{proposition}

The dependence of the lengths  $\vert\mathbf h_1\vert$, $\vert\mathbf h_2\vert$, $\vert\mathbf h_3\vert$ on $\vert\mathbf h_0\vert$ and $\vert\mathbf h_4\vert$ can be explicitly expressed  but it is quite technical, see Remark~\ref{length}.
Before \ls{providing} the proofs of the three previous propositions, we present an illustrative example. 

\begin{example}\mbox{}\rm \label{ex1}
Let us consider two spherical control points 
$$\s_0 = (1,0,0), \qquad   \s_4 = (-0.4330, 0.7500, 0.5000)$$ and the angle $\gamma = 0.6425 \pi$ between them. Let us also choose the spherical control point $$\s_2 = (0.2662,0.8325,-0.4858)$$ which lies on the bisecting great circle $\BC(\s_0,\s_4)$. According to Proposition \ref{pr3} the length of the hodograph control points $\vert \h_0 \vert$ and $\vert \h_4 \vert$ can be arbitrarily chosen. We then set both of them equal to $1$ considering
$$\h_0 = (1,0,0), \qquad   \h_4 = (-0.4330, 0.7500, 0.5000).$$
%Fixing the degree of freedom (see Remark~\ref{fibres}),  
We can then compute the pre-image control points $\A_0$, $\A_2$ as
$$\A_0 = (0,1,0,0), \qquad \A_2 = (-0.4784,0.2338,0.7311,-0.4266)$$ 
so that the first and last equation of \eqref{control_point_rrmf} hold and $\A_0 \, \star \, \A_2$ is a positive multiple of $\s_2$.
The middle equation of \eqref{control_point_rrmf} in turn provides the control point 
$$\h_2 = (0.2338, 0.7311,-0.4266),$$ 
which is a positive multiple of $\s_2$ with length $\vert \h_2 \vert = 0.8782$. We can continue by choosing any point $\s_1\in \BC(\s_0,\s_2)$, as, for example,
$$\s_1 = (0.7686,0.3749,-0.5184).$$ 
The remaining control point of the pre-image $\A_1$ is now uniquely determined as 
$$\A_1 = (-0.3016,0.6819,0.3326,-0.4600)$$ and \eqref{control_point_rrmf} provides $$\h_1 = (0.6819,0.3326,-0.4600), \qquad \h_3 = (-0.1357,0.9250,-0.0188)$$ 
with $\vert \h_1 \vert = 0.8872$, $\vert \h_3 \vert = 0.9351$ and, {consequently}, $\s_3=(-0.1451,0.9892,0.0201).$ %Note that the position of $\s_3$ is uniquely determined by $\s_0$, $\s_1$, $\s_2$, $\s_4$. \cg{[REMOVE?]}

The hodograph $\h(t)$ and the RRMF5-I curve $\r(t)$ for a given $\r(0)$, can be then computed from \ls{\eqref{preimage}--\eqref{PHcontrolpoints}}. The tangent indicatrix $\mathbf t = \mathbf t(t)$ is a rational curve computed via \eqref{rtind}, see Figure~\ref{fig:examples}.

Let us now consider the same spherical control points $\s_0$, $\s_1$, $\s_2$, $\s_4$ and let us define
$$\tilde \h_0 = \h_0 = (1,0,0), \qquad   \tilde \h_4 = (-0.1429,0.2475,0.1650),$$
keeping 
$\vert \tilde \h_0 \vert = \vert \h_0 \vert$ and modifying  $\vert \tilde \h_4 \vert = 0.33 \,\vert \h_4 \vert$. As a consequence, the new pre-image quaternion coefficients are
$$
 \tilde \A_0 \,=\, (0,1,0,0)\,, \quad  
 \tilde \A_1 \,=\, (-0.2286,0.5168,0.2521,-0.3486)\,, \quad
 \tilde \A_2 \,=\, (-0.2748,0.1343,0.4200,-0.2451)\,,  
$$
and the new inner hodograph control points become
$$ \tilde \h_1 \,=\, (0.5168,0.2521,-0.3486)\,, \quad
    \tilde \h_2 \,=\, (0.1343,0.4200,-0.2451)\,,\quad
    \tilde \h_3 \,=\, (-0.0591,0.4027,-0.0082)\,,  
$$
with the lengths $\vert \tilde \h_1 \vert = 0.6725$, $\vert \tilde \h_2 \vert = 0.5045$, $\vert \tilde \h_3 \vert = 0.4071$. Note that the same spherical control point $\s_3$ is obtained due to $\tilde \h_3/|\tilde \h_3|=\h_3/|\h_3|$. Again, we obtain the hodograph $\tilde \h(\tilde t)=\sum_{i=0}^4\tilde \h_iB(\tilde t)$ the RRMF curve $\tilde \r(\tilde t)$ and the indicatrix $\tilde{\mathbf{t}}(\tilde t)$, as before, see Figure~\ref{fig:examples}. The two tangent indicatrices are related via the rational linear reparameterization, $t = t(\tilde t)$ with $t :[0,1]\rightarrow [0,1]$ defined as follows,
\begin{equation*}
    t=\frac{ 0.7579 \, \tilde t}{1-0.2421\, \tilde t}.
\end{equation*}

\begin{figure}[t!] 
\includegraphics[width=0.33\linewidth]{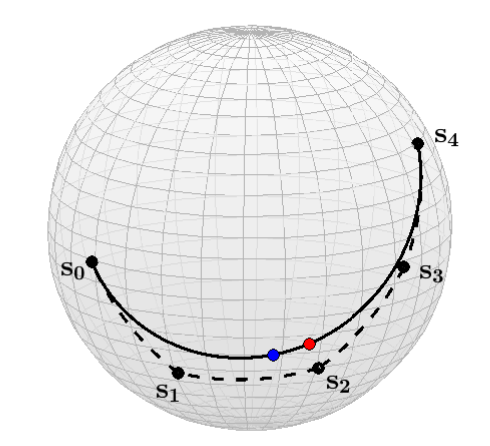}
\includegraphics[width=0.59\linewidth]{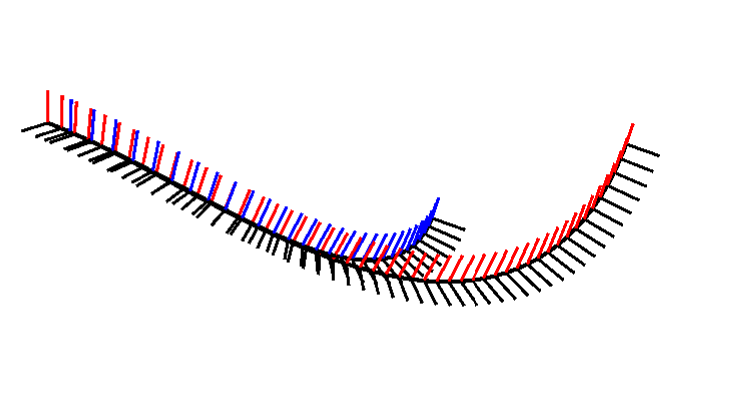}
\caption{Left: the spherical control points $\s_i, i=0,\ldots,4,$ (black dots) of Example~\ref{ex1}, the related spherical control polygon (dashed black line), and the coincident tangent indicatrices $\mathbf{t}(t)$ and $\tilde{\mathbf{t}}(\tilde t\,)$ (black line). The points $\mathbf{t}\left(1/2\right)$ (red dot) and $\tilde{\mathbf{t}} \left(1/2\right)$ (blue dot) are also shown. Right: the corresponding RRMF5-I curves $\r(t)$ and $\tilde \r(\tilde t\,)$ (black lines), sharing the same initial point, with the associated $\f_2$ (black line) and $\f_3$ (red and blue lines for $\r$ and $\tilde \r$, resepectively) 
RMF vectors.}
\label{fig:examples}
\end{figure}   
\end{example}

The symmetrical structure of the hodograph of any RRMF5-I curve, defined by the control points in \eqref{control_point_rrmf}, is highlighted in the following lemma. 
\begin{lemma}\label{lemma}
    Let $\h_b$ and $\h_e$ be two non zero vectors so that $\h_b \times \h_e \neq {\bf 0}$ and consider all the quaternions $\A_b$, $\A_e$ satisfying 
\begin{equation} \label{kv1}
\h_b=\A_b \, \star \, \A_b, \qquad \h_e=\A_e \, \star \, \A_e. \end{equation}
The locus of all vectors 
 \begin{equation}\label{HM}
     \h_m= \A_b \, \star \, \A_e
\end{equation}
is an ellipse lying in the plane identified by $\bis(\h_b,\h_e)$ and $\cross(\h_b,\h_e)$. Moreover, the major axis of the ellipse has the direction of the bisector $\bis(\h_b, \h_e)$ and length $\sqrt{\vert\h_b\vert\vert\h_e\vert}$.
The minor axis, instead, has the perpendicular direction $\cross(\h_b, \h_e)$ and  length $\sin{\frac{\gamma}{2}}\sqrt{\vert\h_b\vert\vert\h_e\vert}$,
where $\gamma$ is the angular distance between $\h_b$ and $\h_e$. 
Consequently, the set of all points $\mathbf h_m$ can be parameterized as
\begin{equation}\label{canPar}
   \mathbf h_m(\varphi)=\sqrt{\vert\h_b\vert\vert\h_e\vert}  \left (    \bis(\h_b, \h_e)\cos \varphi    +\sin{\frac{\gamma}{2}}\,\, \cross(\h_b, \h_e) \sin \varphi \right ), \qquad \varphi \in [0,2\pi).
\end{equation}

\end{lemma}

\begin{figure}[t!]
\begin{center}
    \includegraphics[width=0.40\linewidth]{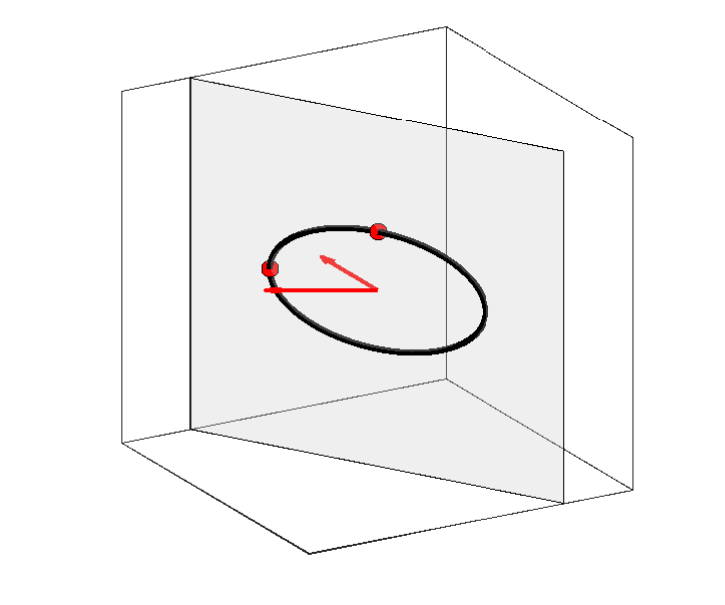}
    \includegraphics[width=0.40\linewidth]{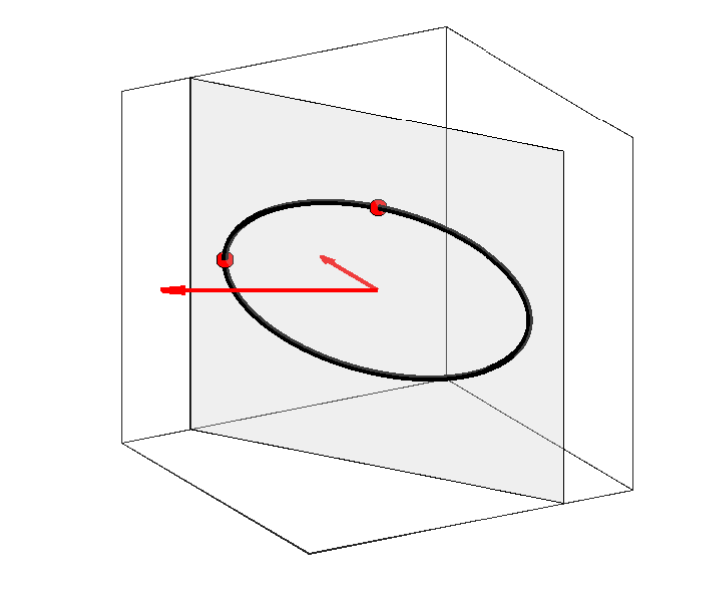}
    \caption{Two vectors $\h_b$ and $\h_e$ (red arrows) with $\vert  \h_b \vert = 1$\, (left) or $\vert  \h_b \vert = 2$ (right), and $\vert  \h_e \vert = 1$. The ellipse that corresponds to the locus of vectors $\h_m$ (black line), as stated in Lemma~\ref{lemma}, and lies on the plane (gray rectangle) identified by $\bis(\h_b,\h_e)$ and $\cross(\h_b,\h_e)$ is also shown. The red dots represent the major and minor axes of the ellipse.}
    \label{fig_ellipses}    
\end{center}
\end{figure}

\begin{proof}\mbox{}
If we define $\hat{\A_b}:=\sqrt{\h_b}\,\,\bis(\ii, \h_b)$, $\hat{\A_e}:=\sqrt{\h_e}\,\, \bis(\ii, \h_e),$ then all solutions of \eqref{kv1}
have the form
    \begin{equation} \label{canonicformA}
        \A_b = \hat{\A_b} e^{\ii \alpha_b} 
        \quad \text{and}\quad  
        \A_e = \hat{\A_e} e^{\ii \alpha_e},
    \end{equation}
    as shown in \eqref{sqrt_quaternion}. Using \eqref{com_mult} and setting $\theta := \alpha_e - \alpha_b$, we can expand \eqref{HM} as
    \begin{equation} \label{hm}
        \h_m =\frac{1}{2}(\hat \A_b e^{\ii \alpha_b} \ii e^{-\ii \alpha_e} \hat \A_e^* + \hat \A_e e^{\ii \alpha_e} \ii e^{-\ii \alpha_b} \hat \A_b^*)=   \hat \A_b \, \star \, \hat \A_e \,\cos \theta+ \hat \A_b \, \square \, \hat \A_e \, \sin\theta,
    \end{equation} 
    which is a parameterization of an ellipse (see Figure~\ref{fig_ellipses}) with conjugate (not necessarily perpendicular) diameters $\hat \A_b \, \star \, \hat \A_e$ and $\hat \A_b \, \square \, \hat \A_e$.
In order to prove that the ellipse lies in the bisecting plane it is sufficient to observe that \cite{jssz2013}
\begin{equation*} 
(\s_b - \s_e) \cdot \hat \A_b \,\square \,\hat \A_e = 0
\quad
\text{and} 
\quad
(\s_b - \s_e) \cdot \hat  \A_b \, \star \, \hat \A_e = 0\, ,
\end{equation*}
where, as before, $\s_b:= \h_b/\vert \h_b \vert$ and  $\s_e:= \h_e/\vert \h_e \vert$. 

We now proceed with the analysis of the axes of the ellipse.  Without loss of generality, we may assume a special position $\ii = \s_b$ obtaining the simplified expressions
 \[\hat \A_b = \sqrt{\vert \h_b \vert }\s_b \qquad \text{and} \qquad \hat \A_e = \sqrt{\vert \h_e \vert }\frac{\s_b + \s_e}{\vert \s_b + \s_e \vert },\]
which imply
 \[ \hat \A_b \, \star \, \hat \A_e = \text{vect}(\hat \A_b \,\s_b\, \hat \A_e^*) = \sqrt{\vert \h_b \vert  \vert \h_e \vert }\frac{\s_b + \s_e}{\vert \s_b + \s_e \vert}=\sqrt{\vert \h_b \vert  \vert \h_e \vert }\bis(\mathbf h_b, \mathbf h_e)
 \]
 and 
  \[ \hat \A_b \, \square \, \hat \A_e = \text{vect}(\hat \A_b \hat \A_e^*) = -\sqrt{\vert \h_b \vert  \vert \h_e \vert }\frac{\s_b \times \s_e}{\vert \s_b + \s_e \vert}=\sqrt{\vert \h_b \vert  \vert \h_e \vert }\sin{\frac{\gamma}{2}}\,\cross(\mathbf h_b, \mathbf h_e),
 \]
where the trigonometric equality
$$\sin{\frac{\gamma}{2}}=\frac{\vert \s_b \times \s_e \vert}{\vert \s_b + \s_e \vert}$$
was considered.
 It is then clear that $\hat \A_b \, \star \, \hat \A_e$ and $\hat \A_b \, \square \, \hat \A_e$, directed as $\bis(\h_b,\h_e)$ and $\cross(\h_b,\h_e)$, respectively, are the perpendicular axes of the ellipse (see Figure~\ref{fig_ellipses}). Setting $\theta=\varphi$ we obtain \eqref{canPar}. 

\end{proof}

\begin{remark} \label{rem2}
Concerning Lemma~\ref{lemma}, if $\h_b$ and $\h_e$ in \eqref{kv1} are replaced by $\h_b^r :=  R_\psi^\ii(\h_b)$ and $\h_e^r := R_\psi^\ii(\h_e)$, with $ R_\psi^\ii$ denoting any rotation about the $\ii$ axis, \ls{some} quaternion computations can be used to show that $\h_m^r = R_\psi^\ii(\h_m)$, where $\h_m$ and $\h_m^r$ are defined as in \eqref{canPar}, respectively referring to the pair $\h_b$, $\h_e$ and $\h_b^r$, $\h_e^r$. This proves that the parameterization \eqref{hm} is covariant with respect to any rotation $R_\psi^\ii$. However, it can also be verified that this feature is not true for a general rotation. Furthermore, the special position $\mathbf s_b=\ii$ adopted in the proof of the lemma can be interpreted as follows. Let $\hat {R}$ be a specific rotation satisfying $\hat{R}(\mathbf s_b)=\ii$. Then, the rotated points $\hat{R}(\mathbf h_b)$, $\hat{R}(\mathbf h_e)$ can be used for the computation of $\hat \A^r_b \, \star \, \hat \A^r_e$ and $\hat A^r_b \, \square \, \hat \A^r_e$ which will be mutually perpendicular, as shown in the proof. Rotating these vectors back using $ \hat R^{-1}$, we obtain the canonical orthogonal parameterization \eqref{canPar}, which does not depend on the particular choice \ls{of} $\hat R$. In fact, any other rotation $R$ so that $R(\mathbf s_b)=\ii$ can be represented as $ R =  R_\psi^\ii \hat{ R }$. On the other hand, in a general position of $\mathbf h_b$, $\mathbf h_e$, the expression \eqref{hm} does not directly provide an orthogonal parameterization since $\hat \A_b \, \star \, \hat \A_e$ and $\hat \A_b \, \square \, \hat \A_e$ are not perpendicular. Note that the relation to the canonical parameterization~\eqref{canPar} is very simple and corresponds to the shift in the parameter values $\varphi=(\theta-\tilde \theta)$ for some fixed angle $\tilde \theta$, see Remark~\ref{bijection} in the next section.
\end{remark}

\vspace{0.001cm}
\subsection{Proofs of the geomertic charaterization and related consequences} \label{sec::proof} 
\subsubsection*{Proof of Proposition~\ref{pr1}}
Lemma \ref{lemma} states that the point $\mathbf h_m$ lies in the bisecting plane of the points $\mathbf h_b$ and $\mathbf h_e$. 
Due to equation \eqref{control_point_rrmf} we can apply Lemma \ref{lemma} three times:
\begin{equation}\label{3cases}
\mathbf h_b=\mathbf h_0, \mathbf h_e=\mathbf h_4, \mathbf h_m=\mathbf h_2, \qquad
\mathbf h_b=\mathbf h_0, \mathbf h_e=\mathbf h_2, \mathbf h_m=\mathbf h_1, \qquad
\mathbf h_b=\mathbf h_4, \mathbf h_e=\mathbf h_2, \mathbf h_m=\mathbf h_3
\end{equation}
and \ls{consider} the normalization $\mathbf s_i=\mathbf h_i/\vert\mathbf h_i\vert$. \ls{This proves} Proposition \ref{pr1}.

\subsubsection*{Proof of Proposition~\ref{pr2}}
In order to decide whether some configuration of the different $\mathbf h_i$ or of their spherical counterparts $\mathbf s_i$ corresponds to an RRMF quintic, we must construct $\mathcal A_0$, $\mathcal A_1$, $\mathcal A_2$ so that \eqref{rrmfstar} and \eqref{control_point_rrmf} are satisfied. Consider three spherical points $\s_0$, $\s_4$ and $\s_2\in \BC(\s_0,\s_4)$ \ls{from} Proposition \ref{pr2} and any corresponding hodograph control points $$\mathbf h_0=\vert\mathbf h_0\vert\s_0,\qquad \mathbf h_4=\vert\mathbf h_4\vert\s_4.$$
Applying Lemma \ref{lemma} to the first case of \eqref{3cases}, any corresponding $\mathbf h_2$ must be of the form 
\begin{equation}\label{h2}
\mathbf h_2=\sqrt{\vert\h_0\vert\vert\h_4\vert}  \left (    \bis(\h_0, \h_4)\cos \varphi_2    +\sin{\frac{\gamma}{2}}\,\, \cross(\h_0, \h_4) \sin \varphi_2 \right ), \qquad \varphi_2 \in [0,2\pi),
\end{equation}
where $\gamma$ is the angle between $\mathbf h_0$ and $\mathbf h_4$.
Only for one specific choice $\varphi_2=\hat \varphi_2$ this point will be a positive multiple of $\s_2$ so that
\begin{equation}\label{h2length}
    \vert\mathbf h_2\vert=\sqrt{\vert\h_0\vert\vert\h_4\vert}\sqrt{(\cos \hat \varphi_2)^2+\left(\sin \frac{\gamma}{2}\right)^2(\sin \hat \varphi_2)^2}
    \end{equation}
and the point $\mathbf h_2=\vert\mathbf h_2\vert\,\s_2$ is uniquely determined. Once the preimage degree of freedom is fixed (see Remark \ref{fibres}), the quaternion coefficients
$$\mathcal A_0=\hat {\mathcal A}_0, \qquad \mathcal A_2=\hat{\mathcal A}_2 e^{\ii \hat \varphi_2}$$ can be determined. The remaining preimage coefficient $\mathcal A_1$ must satisfy \eqref{rrmfstar} and the solutions to this problem, see \eqref{sqrt_quaternion}, form the one parameter system 
$$\mathcal A_1=\sqrt{\vert\mathbf h_2\vert}\,\bis(\ii,\mathbf h_2)\,e^{\ii\theta_1}=:\hat {\mathcal A}_1 e^{\ii\theta_1}$$ 
giving infinitely many pairs via \eqref{control_point_rrmf} 
$$
\mathbf h_1=\A_0\,\star\,\A_1, 
\qquad 
\mathbf h_3=\A_1\,\star\,\A_2,
$$ 
and, consequently, infinitely many pairs $\s_1$, $\s_3$. Moreover, the conditions of Lemma \ref{lemma} are satisfied for the second and third cases of \eqref{3cases} and \eqref{hm}, which yield
    \begin{eqnarray}\label{h1}
      \mathbf h_1&=&  \A_0 \, \star \, \hat \A_1 \,\cos \theta_1+  \A_0 \, \square \, \hat \A_1 \, \sin\theta_1\\ \label{h3}
         \mathbf h_3&=&  \A_2 \, \star \, \hat \A_1 \,\cos \theta_1+ \A_2 \, \square \, \hat \A_1 \, \sin\theta_1.
    \end{eqnarray}
This provides the smooth and orientation preserving bijection $\s_1\leftrightarrow \s_3$, which concludes the proof of Proposition \ref{pr2}. 

\subsubsection*{Proof of Proposition~\ref{pr3}}
We have already seen that the position of $\s_0$, $\s_4$ and $\s_2$ together with  $\vert\mathbf h_0\vert$ and $\vert\mathbf h_4\vert$ fully determine  $\mathcal A_0$ and $\mathcal A_2$. In a similar way, adding the position $\s_1$ (and consequently $\s_3$) the preimage control point $\mathcal A_1$ is uniquely determined. The preimage control points determine the hodograph control points via \eqref{control_point_rrmf}. Note that, for a fixed configuration of the spherical points $\s_i$, if two sets of lengths $\vert \h_0 \vert$, $\vert \h_4 \vert$  and $\vert \tilde \h_0 \vert$, $\vert \tilde \h_4 \vert$ are given, the hodograph control points are related as  $\tilde \h_0 = \mu^2  \h_0$ and $\tilde \h_4 = \mu^2 \lambda^4 \h_4$, with 
\[
\mu = \sqrt{\frac{\vert \tilde \h_0 \vert}{\vert \h_0 \vert}}, 
\qquad 
\lambda = \sqrt[4]{\frac{\vert \tilde \h_4 \vert \vert \h_0 \vert}{\vert \tilde \h_0 \vert \vert \h_4 \vert}}.
\]
We then have
$$ \tilde \A_0 = \mu \A_0, \qquad \tilde \A_1 = \mu \lambda \A_1, \qquad \tilde \A_2 = \mu \lambda^2 \A_2, $$
as specified in the assumption of 
Proposition \ref{PropRepar}. Consequently,  the two tangent indicatrices are related via a linear rational reparameterization and the proof of \ls{Proposition}~\ref{pr3} is concluded.
%\qed

\begin{remark}\label{bijection}
Since the relation  $\s_1\leftrightarrow \s_3$ does not depend on the lengths $\vert \h_0 \vert$, $\vert \h_4 \vert$, we can set $\vert \h_0 \vert=\vert \h_4 \vert=1$. In this case, \eqref{h2length} \ls{implies}
\begin{equation}
    \vert\mathbf h_2\vert=\sqrt{(\cos \varphi_2)^2+ \left(\sin \frac{\gamma}{2} \right)^2(\sin \varphi_2)^2}.
\end{equation}
Moreover, applying again Lemma \ref{lemma}, the two ellipses containing $\mathbf h_1$ and $\mathbf h_3$ will be isomorphic with the lengths of the major and minor axes equal to $\sqrt{\vert\mathbf h_2\vert}$ and $\sin{(\delta/2)}\sqrt{\vert\mathbf h_2\vert}$, respectively, where $\delta$ is the angle between the vectors $\s_0$ and $\s_2$ (the same as the angle between the vectors $\s_4$ and $\s_2$). The parameterizations \eqref{h1} and \eqref{h3}  of these two ellipses do however not preserve this isometry but there is a parametric shift as explained in Remark \ref{rem2}. Choosing the special position $\s_0=\mathbf \ii$ without lost of generality, there will be no shift for $\mathbf h_1$ and the parameterization \eqref{h1} will be canonical (orthogonal). But the shift of \eqref{h3} will remain and we were able to express it using the \ls{Mathematica} software. Indeed, \eqref{h3} is related to the canonical parameterization \eqref{canPar} via $\varphi_1=(\theta_1-\tilde \theta_1)$ where
{\begin{equation}\label{shift}
      \tilde \theta_1= 
   \frac{1}{2}\,  \operatorname{\ls{atan2}} (x_\theta,y_\theta)\,,
\end{equation}}
with 
\begin{equation*}
x_{\theta} = 4\sin \varphi_2 \cos\frac{\gamma}{2}\left(\sin \frac{\gamma}{2} \right)^2 \sqrt{3-\cos \gamma +(1 + \cos \gamma)\cos(2\varphi_2)}\,, \qquad y_\theta =\cos(2\varphi_2)(\sin\gamma)^2 + 4\left(\sin \frac{\gamma}{2} \right)^4\,.
\end{equation*}
\end{remark}   

\begin{remark}\label{length} 
For a given configuration of spherical control points $\s_i$, the length $\vert\mathbf h_2\vert$ is determined by $\vert\mathbf h_0\vert$ and $\vert\mathbf h_4\vert$ via \eqref{h2length}. Considering again the special position $\s_0=\ii$, the parameterization \eqref{h1} becomes canonical (see Remark~\ref{bijection}) and we can then apply the formula analogue to \eqref{h2length} to get

\begin{eqnarray*}
    \vert\mathbf h_1\vert&=&\sqrt{\vert\h_0\vert\vert\h_2\vert}\sqrt{(\cos \theta_1)^2+(\sin \frac{\delta}{2})^2(\sin \theta_1)^2}\\
    &=&\sqrt[4]{\vert\h_0\vert^3\vert\h_4\vert} \sqrt[4]{(\cos \varphi_2)^2+(\sin \frac{\gamma}{2})^2(\sin \varphi_2)^2}\sqrt{(\cos \theta_1)^2+(\sin \frac{\delta}{2})^2(\sin \theta_1)^2},
\end{eqnarray*}
where $\delta$ is the angle between the vectors $\s_0$ and $\s_2$.

The length $\vert\mathbf h_3\vert$ is given by a very similar formula but the shift \eqref{shift} must be \ls{considered}
\begin{eqnarray*}
    \vert\mathbf h_3\vert&=&\sqrt{\vert\h_2\vert\vert\h_4\vert}\sqrt{\cos (\theta_1-\tilde \theta_1)^2+(\sin \frac{\delta}{2})^2\sin (\theta_1-\tilde \theta_1)^2}\\
    &=&\sqrt[4]{\vert\h_0\vert\vert\h_4\vert^3} \sqrt[4]{(\cos \varphi_2)^2+(\sin \frac{\gamma}{2})^2(\sin \varphi_2)^2}\sqrt{\cos (\theta_1-\tilde \theta_1)^2+(\sin \frac{\delta}{2})^2\sin (\theta_1-\tilde \theta_1)^2}.
\end{eqnarray*}

\end{remark}

\section{$G^1$ Hermite interpolation with one sided frame conditions}\label{sec:local_alg}

On the basis of \ls{previous} geometric considerations, we now introduce an algorithm for the construction of an RRMF5-I curve interpolating  suitable $G^1$ Hermite data at the end points, together with an initial frame orientation. We then need to take into account the interpolation of an initial and final points, $\p_i$ and $\p_f$, an initial frame ($\u_i$, $\mathbf{v}_i$, $\w_i$ ) together with the final tangent direction $\u_f$. Note that we consider a specific configuration where the non vanishing displacement $\mathbf{\Delta} \p:= \p_f - \p_i$ and the corresponding projection on the unit sphere $\mathbf{\Delta} \u := \mathbf{\Delta} \p\, /\,\vert \mathbf{\Delta} \p \vert$ satisfy
\begin{equation} \label{construction_hp}
\u_i \cdot \mathbf{\Delta} \u = \mathbf{\Delta} \u \cdot \u_f\,.
\end{equation}
\ls{This condition ensures that the unit dispacement $\mathbf{\Delta} \u$ lies on $\BC(\s_0,\s_4)$, that will be later suitably exploited in the development of our algorithm, as detailed below, see Subsection~\ref{PHdisp}}.

\subsection{Local algorithm} \label{alg1}
By considering the general PH parametric form provided by \eqref{PHcurve} and \eqref{PHcontrolpoints}, we want to construct an RRMF5-I curve $\r$ that satisfies the end-point interpolation conditions
\begin{equation} \label{C0interp}
\mathbf{r}(0) = \p_i, 
\quad
\mathbf{r}(1) = \p_f, 
\end{equation}
together with the initial frame orientation
\begin{equation} \label{frameInterp}
\qquad (\f_1(0),\f_2(0),\f_3(0)) = (\u_i, \mathbf{v}_i,\w_i),  
\end{equation}
and the $G^1$ interpolation conditions
\begin{equation} \label{G1interp}
\mathbf{r}'(0)
= \mu^2\, \u_i,\qquad
\mathbf{r}'(1)
= \mu^2 \, \u_f ,
\end{equation}
for any real positive $\mu$, with $\u_i \times \u_f \neq {\bf 0}$. In theory two different values $\mu_i$ and $\mu_f$ could be used in \eqref{G1interp} instead of $\mu$. Our choice simplifies the computations while simultaneously being compatible with arc length oriented parameterizations. The interpolation of the assigned end tangent directions  corresponds to the following two vector conditions:
 \begin{equation*} 
  \A_0\,\ii\,\A_0^* = \mu^2 \u_i\, , \qquad
 \A_2\,\ii\,\A_2^* =\mu^2 \u_f\, .
  \end{equation*}
This leads to 
\[ 
\A_0 = \mu \, \U_0 \qquad \text{and} \qquad \A_2 = \mu \, \U_2\ls{\,,}
  \]
where 
\begin{equation} \label{U0U2def}
\U_0 = \b(\u_i,\ii) \,e^{\alpha_0 \ii},
\qquad 
\U_2 = \b(\u_f,\ii)\,e^{\alpha_2 \ii} = \b(\u_f,\ii)\, e^{\alpha_0 \ii}\,e^{\theta_2\ii} =: \hat \U_2 \,e^{\theta_2 \ii},
\end{equation}
$\theta_2 := \alpha_2 - \alpha_0$, and we assume $\ii \neq -\u_i,\,\ii \neq -\u_f$. As shown in \cite{fgms2012}, the interpolation of the initial frame orientation can be obtained if the angular parameter $\alpha_0$ satisfies
\begin{align} \label{alfa0}
&\cos(2\alpha_0) = \k_0 \cdot \w_i \qquad  \text{with} \qquad \k_0 = 2(\k \cdot \b(\u_i,\ii))\b(\u_i,\ii) - \k, \nonumber \\
&\sin(2\alpha_0) = -\j_0 \cdot \w_i \qquad  \text{with} \qquad \j_0 = 2(\j \cdot \bis(\u_i,\ii))\b(\u_i,\ii) - \j.
\end{align}
Note that the first and last spherical control points of the hodograph satisfy 
\begin{equation}\label{t0t4}
    \s_0 = \U_0 \,\ii \,\U_0^* = \u_i, 
    \quad \s_4 = \U_2 \,\ii \,\U_2^* = \u_f
\end{equation}
and, in view of the hypothesis \eqref{construction_hp}, the spherical displacement $\mathbf{\Delta} \u$ lies on the great circle $\BC(\s_0,\s_4)$. Moreover,  we know from Proposition \ref{pr1} that also the spherical control point $\s_2$ of the hodograph of any RRMF5-I lies on the great circle $\BC(\s_0,\s_4)$ and, according to \eqref{hm}, its  position only depends on the free angular parameter $\theta_2$, resulting in 
\begin{equation}\label{s2}
    \s_2(\theta_2) = \frac{\U_0 \star \U_2(\theta_2)}  {\vert \U_0 \star \U_2(\theta_2) \vert}\,.
\end{equation} 
From the RRMF condition \eqref{rrmf}, we can then compute  
\begin{equation*}
\mathcal A_1(\theta_2)= \mu \sqrt{\vert \U_0 \, \star \, \U_2(\theta_2) \vert} \, \U_1(\theta_2), \quad \text{with} \quad \U_1(\theta_2) = \bis(\ii, \s_2(\theta_2))\,e^{\alpha_1 \ii}=\bis(\ii, \s_2(\theta_2))\,e^{\alpha_0 \ii}\,e^{\theta_1 \ii}=: \hat \U_1(\theta_2) \,e^{\theta_1 \ii}\,,
\end{equation*} 
where only the dependence on the free angular parameter $\theta_2$ is highlighted, since $\theta_1:= \alpha_1 - \alpha_0$ is fixed as detailed below.  
In view of Proposition \ref{pr2}, this angular value  fixes the position of $\s_1$ and $\s_3$ on $\BC(\s_0,\s_2)$ and $\BC(\s_4,\s_2)$, respectively. We select its value so that
\begin{equation}\label{alpha1} 
    \text{scal}\left( (\U_0 + \U_2) \,\ii \,\U_1^* \right) = 0\,.
\end{equation}
The above condition can be interpreted as follows. We require that the rotation angle of the two unit quaternions $\U_1 \,\ii \, \U_0^*$ and $\U_2 \,\ii \, \U_1^*$, which rotate $\s_0$ and $\s_2$ into $\s_2$ and $\s_4$, respectively, is the same. Another interesting consequence of this choice for $\theta_1$ can be underlined. Let us consider the unit quaternion
\begin{equation}\label{sm}
    \s_m :=\frac{(\U_0+\U_2)\,\ii \,\U_1^* }{ \vert \U_0 + \U_2 \vert},
\end{equation}
which rotates $\s_2$ on the unit sphere into 
\begin{equation}\label{s02}
    \s_{02} := \frac{(\U_0+\U_2)\,\ii \,(\U_0+\U_2)^* }{ \vert \U_0 + \U_2 \vert^2} = \frac{\s_0 + \s_4 + 2 \ \U_0 \star \U_2 }{\vert \U_0 + \U_2 \vert^2} = \frac{\vert \s_0 + \s_4\vert \, \bis(\s_0,\s_4) + 2 \ \vert \U_0 \star \U_2\vert \s_2}{\vert \U_0 + \U_2 \vert^2} .    
\end{equation}
We can observe that $\s_{02} \in \BC(\s_0,\s_4),$ since both $\s_2$ and  $\bis(\s_0,\s_4)$ belong to $\BC(\s_0,\s_4)$. 
Taking into account that the scalar part of a unit quaternion can be seen as the cosine of the semi-angle of rotation, the previous equation points out that $\s_m$ is a pure vector quaternion, see (\ref{alpha1}), and identifies a rotation of an angle $\pi$ between $\s_2$ and $\s_{02}$, namely $\s_m = \pm \bis(\s_{02},\s_2)$. A direct consequence is that 
\begin{equation}\label{tm}
    \s_m \in \BC(\s_0,\s_4)\,.
\end{equation}
Solving equation \eqref{alpha1} it is possible to find two values for the parameter $\theta_1$, given by
\begin{equation}\label{theta1}
\operatorname{atan2}(x_1 , y_1) \quad\text{and}\quad \operatorname{atan2}(x_1 , y_1) + \pi\,,
\end{equation}
with 
\begin{equation}\label{xyfortheta}
x_1 = \text{scal}\left( \U_0 \,\ii \, \hat \U_1^* - \hat \U_1 \,\ii \, \U_2^*\right)\,, \qquad y_1 = -\,\text{scal}\left(\U_0 \, \hat \U_1 ^* + \hat \U_1\,\U_2^* \right),
\end{equation}
corresponding to the two admissible directions for  $\s_m$. To obtain a better shape of the spherical control polygon identified by $\s_0,\ldots,\s_4,$ we choose the solution for $\theta_1$ so that 
\begin{equation} \label{smdef}
\s_m = \bis(\s_{02},\s_2)\,.
\end{equation}
Finally, with the remaining free parameters, the angle $\theta_2$ and the positive coefficient $\mu,$  we want to \ls{fulfill} the end point condition
\begin{equation}\label{alg_integral}
 \mathbf{\Delta} \p = \int_0^1 \A(t) \ii A(t)^* dt = \mu^2  \int_0^1 \mathcal{V}(t) \ii \mathcal{V}(t)^* dt =: \frac{\mu^2}{5}\, \mathbf{I}_\gamma(\theta_2) ,  
\end{equation}
with
\[
  \mathcal{V}(t) := \U_0 (1-t)^2 + \sqrt{\vert \U_0 \, \star \, \U_2 \vert}\, \U_1 \,2t(1-t) + \U_2t^2,
\]
recalling that $\gamma$ is the angular distance between $\u_i$ and $\u_f.$ Assuming that $\mathbf{I}_\gamma(\theta_2)$ does not vanish, the positive parameter $\mu$ can be used to ensure that the two vector quantities on the two sides of (\ref{alg_integral}) have the same module, i.e.,
\begin{equation} \label{lamval}
\mu = \sqrt{\frac{5 \vert \mathbf{\Delta} \p \vert }{\vert \mathbf{I}_\gamma(\theta_2) \vert}}.
\end{equation}
Consequently, the existence of solutions depends on the possibility to select the free angle $\theta_2$ so that $\mathbf{I}_\gamma(\theta_2)$ is a nonzero vector aligned with $\Delta \p.$
In the next subsection we then present an analysis on the variation  of the direction associated to $\mathbf{I}_\gamma(\theta_2)$ when $\theta_2$ varies, and we denote it as the \emph{scaled PH displacement}. This study is useful to clarify for which values of $\gamma$ the considered interpolation problem admits exactly one solution.  It also motivates the first order Hermite data definition that will be presented in Section~\ref{sec:global_alg} to guarantee the existence of a suitable  $G^1$  RRMF5-I spline curve interpolating an input point stream.

\subsection{Study of the scaled PH displacement} \label{PHdisp}
Recalling equations \eqref{control_point_rrmf} and \eqref{sm}, we can explicitly compute the scaled PH displacement as
 \begin{equation} \label{intvdef}
 \mathbf{I}_\gamma(\theta_2) = \, \q_1 + \q_2 + \q_3 
 \end{equation}
with
\begin{equation} \label{vdef}
\q_1 := \u_i +  \u_f \,, \quad \q_2 = \q_2(\theta_2) :=  \U_0 \star \U_2\,,  \quad
  \q_3 = \q_3(\theta_2) := \sqrt{\vert \U_0 \, \star \, \U_2 \vert}\, (\U_0 + \U_2) \star \U_1\,.
\end{equation}
Among these three vectors, $\q_1 $ (which can not vanish in view of the assumption $\u_f \ne - \u_i$) is the only one not depending on $\theta_2$. We have
 \begin{equation} \label{q1def}
  \q_1 =\vert \u_i+\u_f \vert \, \b  = \sqrt{2(1+ \cos\gamma)}\, \b\,,
  \end{equation}
 where $\b = \bis(\s_0,\s_4)\,,$ and we assume the angular distance $\gamma := \arccos (\u_i \cdot \u_f)$ between $\u_i$ and $\u_f,$ so that 
 $\gamma \in (0\,,\,\pi).$ 
 Also $\q_2$ can not vanish since $\U_0 \star \U_2$ is the vector part of the unit quaternion $\U_0\,\ii\,\U_2^*$ associated to the rotation which maps $\u_f$ into $\u_i \ne \pm \u_f.$ Without loss of generality, we set $\ii = \u_i$ and $\j = -\mathbf{v}_i$, and, consequently $\k = -\w_i$. From (\ref{U0U2def}) and (\ref{alfa0}) we can then set $\U_0 = \ii$ and 
\begin{equation} \label{U2def}
\U_2 = \cos \theta_2 \, \left( \begin{array}{c} 0 \cr \b \end{array} \right) 
\,+\, \sin \theta_2 \ \left( \begin{array}{c} - \cos \frac{\gamma}{2} \cr  \sin \frac{\gamma}{2} \cross \cr \end{array} \right),
\end{equation}
where $\ls{\cross = \cross(\s_0,\s_4)\,,}$ as already considered in (\ref{ndef}). 
We then obtain
$$\U_0\,\ii \,\U_2^* =   \cos \theta_2 \left( \begin{array}{c} 0 \cr \b \end{array} \right) \,+\, \sin \theta_2 \ \left(\begin{array}{c} \cos\frac{\gamma}{2} \cr  \sin\frac{\gamma} {2}\ \cross \end{array} \right),$$ 
 which implies 
 \begin{equation} \label{q2def}
 \q_2 = \q_2(\varphi_2)  = \cos \varphi_2 \ \b + \sin \varphi_2 \sin\frac{\gamma} {2}\ \cross\,,
 \end{equation} 
 where $\theta_2$ has been replaced \ls{by}  $\varphi_2$ because the above expression is an occurrence of \eqref{canPar}, see Lemma \ref{lemma} in Section 2. 
 Consequently, $\q_2(\varphi_2 \pm \pi) = -\q_2(\varphi_2)$, $\q_2(2\pi - \varphi_2) \cdot \cross = - \q_2(\varphi_2) \cdot \cross$, and  $\vert \q_2 \vert = \sqrt{1 - \sin^2 \varphi_2 \ \cos^2 \frac{\gamma}{2}}.$  Recalling that $\s_2$ is the unit vector aligned with $\q_2$ we also have
$$\s_2 = \s_2(\varphi_2) = \frac{\q_2}{\vert \q_2 \vert} = \frac{\cos\varphi_2\ \b \,+\, \sin \varphi_2\ \sin\frac{\gamma}{2} \ \cross}{\sqrt{1 - \sin^2 \varphi_2 \ \cos^2 \frac{\gamma}{2}}}\,.$$ 
  
Note that we can also write
\begin{equation} \label{v3first}
\q_3 = \q_3(\varphi_2) = \sqrt{\vert \U_0 \, \star \, \U_2 \vert}\,  (\U_0 + \U_2) \ii \U_1^* \,,
\end{equation}
since, due to \eqref{alpha1}, the right hand side defines a pure vector quaternion. This means that also $\q_3$  can not vanish because $\U_2 \ne - \U_0,$ being $\gamma \in (0\,,\,\pi).$ Taking into account (\ref{smdef}), we may write
\begin{equation} \label{q3final}
\q_3= \q_3(\varphi_2) = \sqrt{\vert \U_0 \, \star \, \U_2 \vert}\,\vert \U_0+\U_2 \vert \,\s_m \,=\, \sqrt{\vert \U_0 \, \star \, \U_2 \vert}\,\vert \U_0+\U_2 \vert \, \bis(\s_{02},\s_2)\,,
\end{equation}
with $\s_{02}$ defined in (\ref{s02}). Now, since $\U_0 = \ii$ and $\U_2$ is defined as in (\ref{U2def}), we obtain
$$ \left\vert\, \U_0 + \U_2\,\right\vert^2 = \U_0\,\U_0^* + \U_2\,\U_2^* + 2\, \rm{scal}(\U_0\,\U_2^*) = 2\, \ls{\left(1 + \cos\varphi_2 \cos \frac{\gamma}{2}\right)\,.}$$
Thus (\ref{s02}) can be rewritten as 
\begin{equation} \label{s02def}
\s_{02} =  \s_{02}(\varphi_2) =  
\frac{\q_1 + 2\,\q_2(\varphi_2)}{2\,(1 + \cos\varphi_2  \cos \frac{\gamma}{2} )} \,.
\end{equation}
Consequently, $\s_{02}(2\pi-\varphi_2) \cdot \cross = - \s_{02}(\varphi_2) \cdot \cross$
and (\ref{q3final}) then implies that $\q_3(2\pi-\varphi_2) \cdot \cross  = - \q_3(\varphi_2) \cdot \cross.$ 
Thus we can derive preliminary information related to the locus described on the sphere by the central projection of $\mathbf{I}_\gamma(\varphi_2)$ when $\varphi_2$ varies in $[0\,,\,2\pi\ls{)}.$
Since the three vectors $\q_1\,,\,\q_2\,,\, \q_3$ are \ls{a} linear combination of $\b$ and $\cross,$ 
 \begin{equation} \label{Sdef}
 \S_\gamma(\varphi_2) := \frac{\mathbf{I}_\gamma(\varphi_2) }{ \vert \mathbf{I}_\gamma(\varphi_2) \vert } \in \BC(\s_0,\s_4).
 \end{equation}
The next two propositions present key results concerning $\mathbf{I}_\gamma$.

\begin{proposition}  \label{propB}
 Let $\S_\gamma(\varphi_2) $ denote the unit PH displacement introduced in \eqref{Sdef}, expressed as a function of the angle $\varphi_2.$ For all $\gamma \in (0\,,\,\pi)$, we have $\S_\gamma(0) = \b.$  In addition, if $\gamma > 2\pi/5 \, (\gamma< 2\pi/5)$, then $\S_\gamma( \pi) = - \b \,(\S_\gamma(\pi) =  \b),$ while, if $\gamma = 2\pi/5$ then $\mathbf{I}_\gamma(\pi) = {\bf 0}. $ 
 \end{proposition}
 \begin{proof} Formula (\ref{q2def}) implies that $\q_2 (0) = \s_2(0) = \b.$   Considering the definition in (\ref{s02def}), for $\varphi_2 = 0$ it also holds $\s_{02} = \s_{02}(0) = \b.$ Now, in view of (\ref{q3final}), we know that   $\q_3(0)/ \vert \q_3(0) \vert  = \s_m(0) = \bis(\s_{02}(0), \s_2(0)).$ Thus, considering that both $\s_{02}(0)$ and  $\s_2(0)$ coincide with $\b,$ we can conclude that $\q_3(0)/ \vert \q_3(0) \vert = \b.$  Since  $\mathbf{I}_\gamma(0)$ is the sum of three non vanishing vectors all aligned with $\b,$ the first part of the statement is proved. When $\varphi_2 = \pi$, we have
$\q_2 (\pi) = \s_2(\pi) = -\b.$  Thus, from (\ref{s02def}) and (\ref{q1def}), we obtain
$$\s_{02}(\pi) = \frac{\q_1 - 2\b}{2(1 -  \cos \frac{\gamma}{2} )} = \frac{(\sqrt{2(1+ \cos\gamma)} -2)\ \b}{2(1 -  \cos \frac{\gamma}{2} )}\,.$$
Since $\s_m = \bis(\s_2,\s_{02}),$ we get
$$\s_m(\pi) = \frac{(\sqrt{2(1+ \cos\gamma)} -4 + 2 \cos \frac{\gamma}{2} )\ \b}{4(1 -  \cos \frac{\gamma}{2} )}\,.$$ 
Now, considering that $\vert\, \U_0 \, \star \, \U_2(\pi) \,\vert =1$ and $\vert \,\U_0 + \U_2(\pi) \,\vert = \sqrt{2\,(1- \cos \frac{\gamma}{2})}\,,$ from (\ref{q3final}) we obtain 
$$\q_3(\pi) =  \, \frac{(\sqrt{1+ \cos\gamma} -2\sqrt{2} + \sqrt{2} \cos \frac{\gamma}{2} )\ \b}{2\ \sqrt{1 -  \cos \frac{\gamma}{2}} }\,.$$ 
Using $\cos \gamma = 2 \cos^2  \frac{\gamma}{2} -1\,,$ some easy computations lead to
$$ \q_3(\pi) = - \sqrt{2(1 -  \cos \frac{\gamma}{2})} \, \b\,.$$ 
From (\ref{intvdef}) we then have
\begin{equation} \label{exppi}
\mathbf{I}_\gamma(\pi)  = \left( 2 \cos \frac{\gamma}{2} -1 - \sqrt{2(1 -  \cos \frac{\gamma}{2})} \,\, \right) \b \,. 
\end{equation}
\ls{The solution of the equation $2 \cos \frac{\gamma}{2} -1 - \sqrt{2\left(1 -  \cos \frac{\gamma}{2}\right)} = 0$ for $\gamma \in (0\,,\,\pi)$, is $\gamma = 2\pi/5$.} From the associated inequality, it easily descends also the remaining proof of the last part of the statement. \end{proof} 

\noindent
\ls{Introducing the notation 
\begin{equation} \label{sigmaN}
\rho^{\mathbf{b}}_{\gamma}(\varphi_2) := \mathbf{I}_\gamma(\varphi_2) \cdot \bis\,,
\qquad 
\rho^{\mathbf{n}}_{\gamma}(\varphi_2) := \mathbf{I}_\gamma(\varphi_2) \cdot \cross\,, 
\end{equation} 
with $\bis = \bis(\s_0,\s_4),\, \cross = \cross(\s_0,\s_4)\,,$ respectively defined accordingly with \eqref{bdef} and \eqref{ndef}, in the next proposition we start studying the qualitative behavior of $\rho_\gamma^{\cross}(\varphi_2)$.}

\ls{\begin{proposition} \label{propN} 
For any $\gamma \in (0\,,\,\pi)$ and any $\varphi_2 \in [0\,,\,2\pi)$, the following two inequalities hold true,
$$\rho^\mathbf{n}_{\gamma}(\varphi_2 ) (\q_2(\varphi_2) \cdot \cross) \ge 0\,, \quad \rho^\mathbf{n}_{\gamma}(\varphi_2 + \pi)  \, \rho^\mathbf{n}_{\gamma}(\varphi_2 ) \le 0\,,$$ 
where both of them are strict inequalities if $\q_2(\varphi_2) \cdot \cross \ne 0\,.$  
\end{proposition}}
\begin{proof}
First, let us simplify the notation omitting the superscript for referring to the considered dot product. Equations (\ref{intvdef}) and (\ref{vdef}) imply $\rho^\mathbf{n}_{\gamma}(\varphi_2) = ( \q_2(\varphi_2)  +  \q_3(\varphi_2) ) \cdot \cross$. 
Now, $\q_3$ is aligned with $\s_m = \bis(\s_{02},\s_2)$ and, considering the definition of $\s_{02}$  in (\ref{s02def}), $(\s_{02}(\varphi_2) \cdot \cross) (\q_2(\varphi_2) \cdot \cross) \ge 0.$ We also have
$(\q_3(\varphi_2) \cdot \cross)  (\q_2(\varphi_2) \cdot \cross) \ge 0.$ This means that $\rho^\mathbf{n}_{\gamma}(\varphi_2) (\q_2(\varphi_2) \cdot \cross) \ge 0$ and these inequalities are all strict if $\q_2(\varphi_2) \cdot \cross \ne 0.$ Since $\q_2(\varphi_2 + \pi) \cdot \cross = - \q_2(\varphi_2 ) \cdot \cross\,,$ the proof is completed.    \end{proof} 

Taking into account that $\q_1 \cdot \cross = 0$  and $\q_i(2\pi-\varphi_2) \cdot \cross = - \q_i(\varphi_2) \cdot \cross$ for both $i=2$ and $i=3$ and also that, with some additional computations, $\vert \mathbf{I}_\gamma(2\pi-\varphi_2)\vert = \vert \mathbf{I}_\gamma(\varphi_2)\vert\,,$ we can  preliminary note that $\S_\gamma(2\pi- \varphi_2) \cdot \cross = -\S_\gamma(\varphi_2) \cdot \cross.$ This means that the portion of the great circle $\BC(\s_0,\s_4)$ covered by  $\S_\gamma(\varphi_2)$ when $\varphi_2$ varies in $[0\,,\,2\pi]$ is symmetric with respect to the straight line joining $\b$ to $-\b.$

Let us now focus on the case when $\mathbf{I}_\gamma(\varphi_2) $ vanishes.
%------------
\begin{proposition}
When $\gamma \ne 2\pi/5\,,$ the PH displacement $\mathbf{I}_\gamma (\varphi_2)$ does not vanish, regardless from $\varphi_2.$ When $\gamma = 2\pi/5\,,$ it vanishes only for $\varphi_2 = \pi.$ 
\end{proposition}
\begin{proof}In order to get a vanishing PH displacement, we need to ensure the vanishing of both $\rho^\mathbf{n}_{\gamma}(\varphi_2)$ and $\rho^\b_{\gamma}(\varphi_2)$ defined in \eqref{sigmaN}. Now, in view of the previous proposition, $\rho^\mathbf{n}_{\gamma}(\varphi_2)$ vanishes if and only if $\q_2(\varphi_2) \cdot \cross$ vanishes, i.e., only if  $\varphi_2 = 0, \pi.$ On the other hand, due to Proposition \ref{propB},  $\mathbf{I}_\gamma(0)$ can not vanish. Note that the expression of $\mathbf{I}_\gamma(\pi)$ reported in (\ref{exppi}) can vanish just for $\gamma = 2\pi/5.$ \ls{This concludes the proof.}
\end{proof}

\begin{figure}[t!]
\begin{center}
    \includegraphics[width=0.49\linewidth]{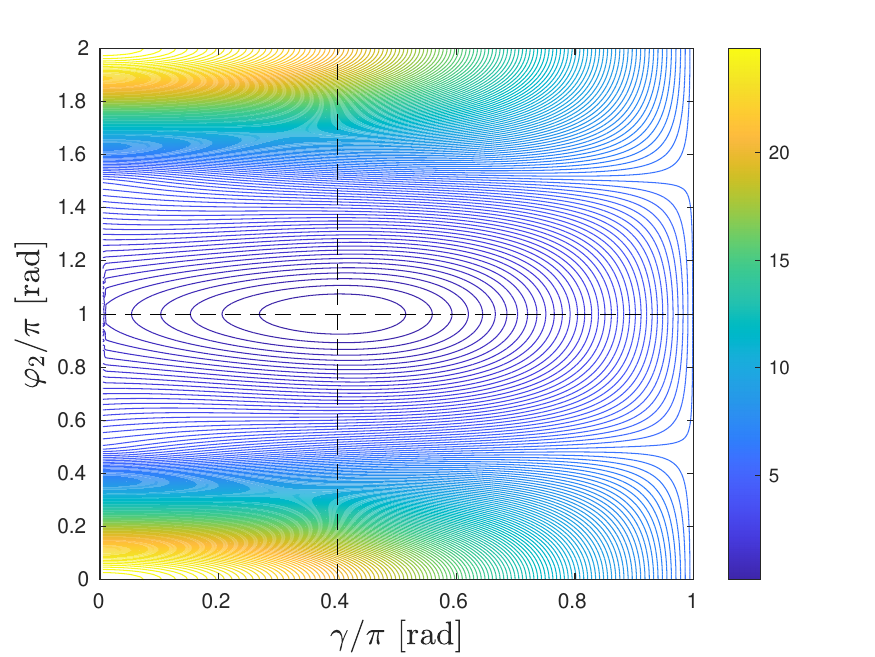}
    \includegraphics[width=0.49\linewidth]{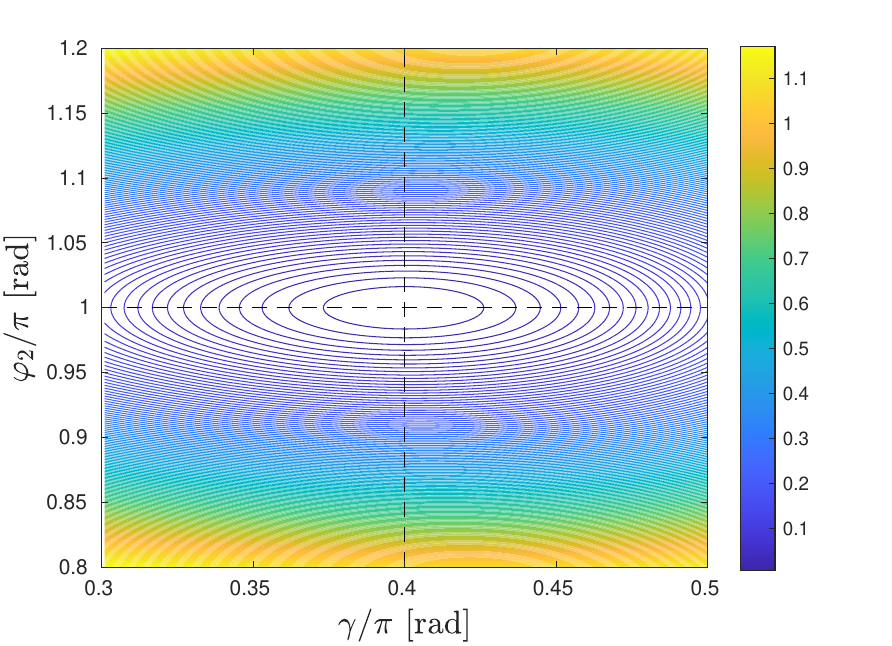}
    \caption{The contour of $\vert \mathbf{I}_\gamma (\varphi_2)\vert$  with $\gamma \in [0 \,,\,\pi]$, $\varphi_2 \in [0 \,,\,2\pi]$ (left) or $\gamma \in [0.3\pi\,,\,0.5\pi]$, $\varphi_2 \in [0.8\pi\,,\, 1.2\pi]$ (right). The intersection of the two dotted lines corresponds to the couple $(\gamma = 2\pi/5 \,,\, \varphi_2 = \pi)$ for which $\mathbf{I}_\gamma (\varphi_2)$ vanishes.}
    \label{integral_contour}  
\end{center}
\end{figure}

\noindent
The above result is also visualized in Figure~\ref{integral_contour}, where the contour lines of $\vert \mathbf{I}_\gamma (\varphi_2) \vert$, considered as a bivariate function of the angle pair $(\gamma,\varphi_2) \in (0\,,\,\pi) \times [0\,,\,2\pi)$, are shown.
%------------------
Unfortunately, $\S_\gamma(\varphi_2)$ does not necessarily cover the full circle $\BC(\s_0,\s_4)$ when $\varphi_2$ varies in $[0\,,\,2\pi),$ regardless of $\gamma.$ The next two propositions are of interest on this concern. 
 \begin{proposition} \label{propF}
 If $\gamma > 2\pi/5,$ then  $\S_\gamma(\varphi_2)$ covers the full circle $\BC(\s_0,\s_4)$  when $\varphi_2$ varies in $[0\,,\,2\pi).$ 
 \end{proposition}
\begin{proof} We have already proved in Proposition \ref{propB} that $\S_\gamma(0)  = \b$ and $\S_\gamma(\pi)  = - \b$ when $\gamma > 2\pi/5\,.$ Thus, for continuity, $\S_\gamma(\varphi_2)  \cdot \b$ assumes all values between  $1$  and $-1$ for $\varphi_2 \in [0\,,\,\pi],$ and analogously in $ [\pi\,,\,2\pi]$. This implies, again for continuity, that there exists $\varphi_{2,1} \in (0\,,\,\pi)$ and  $\varphi_{2,2} \in (\pi\,,\,2\pi)$ such that $ \S_\gamma(\varphi_{2,1}) \cdot \bis =  \S_\gamma(\varphi_{2,2}) \cdot \bis = 0\,.$  Then Proposition \ref{propN} implies that $\rho^\mathbf{n}_{\gamma}(\varphi_{2,1}) = \pm 1$ and $\rho^\mathbf{n}_{\gamma}(\varphi_{2,2}) = \mp 1$ and the proof can be obtained due to continuity arguments. \end{proof}

 \begin{proposition} \label{prop:gamma}
 Let $\BC_h(\s_0,\s_4)$ be the semicircle defined as the half of $\BC(\s_0,\s_4)$ passing  through the spherical points $\cross$, $-\cross$ and $\bis$. If $0 < \gamma \leq 2\pi/5$, then $\S_\gamma(\varphi_2)$ covers a portion of $\BC_h(\s_0,\s_4)$ when $\varphi_2$ varies in $[0\,,\,2\pi).$ 
 \end{proposition}
\begin{proof}  
In order to prove the proposition, it is necessary to verify that the quantity $\rho^\b_{\gamma}(\varphi_2)$ defined in \eqref{sigmaN} is nonnegative if $0 < \gamma \leq 2\pi/5\,.$ then, in view of \eqref{intvdef}, we need to compute the three dot products $\q_1 \cdot \bis,\,\q_2 \cdot \bis, \,\q_3 \cdot \bis\,.$ It can be easily verified that $\q_1 \cdot \b = \vert \q_1 \vert\,, \q_2 \cdot \bis = \cos \varphi_2\,.$ Let consider the definition for $\q_2$ in \eqref{vdef} and compute the bisector $\bis(\s_{02},\s_2)$ via \eqref{bdef}, with $\s_2$ and $\s_{02}$ expressed as in \eqref{s2} and \eqref{s02}, respectively. Substituting them into \eqref{q3final}, after some calculations, we obtain 
\begin{equation} \label{q3b}
\q_3 \cdot \bis = \frac{\vert \q_1 \vert \,\vert \q_2 \vert +\left( \, 2\vert \q_2 \vert +   \vert \U_0 + \U_2 \vert^2 \,\right)\cos \varphi_2}{\sqrt{2\,\vert \q_2 \vert \left( \, \vert \U_0 + \U_2 \vert^2+2\vert \q_2 \vert \, \right)+2\,\vert \q_1 \vert\,\cos \varphi_2}}\,.
\end{equation}
Substituting $\vert \q_1 \vert = \sqrt{2(1+ \cos\gamma)}$, $\vert \q_2 \vert = \sqrt{1 - \sin^2 \varphi_2 \ \cos^2 \frac{\gamma}{2}}$ and $ \vert \U_0+ \U_2\vert^2 = 2 (1 + \cos\varphi_2 \cos \frac{\gamma}{2})$ into \eqref{q3b}, we derive an explicit expression for $\rho^\b_{\gamma}(\varphi_2)$ as 
\begin{footnotesize}
\begin{multline*}
\rho^\b_{\gamma}(\varphi_2) = \sqrt{2(1+ \cos\gamma)} +\cos \varphi_2 +\frac{\sqrt{2 + 2\cos\gamma}\sqrt{1 - \sin^2\varphi_2\cos^2\frac{\gamma}{2}} + \left(2 + 2\cos \varphi_2 \cos \frac{\gamma}{2} + 2\sqrt{1 - \sin^2 \varphi_2 \cos^2 \frac{\gamma}{2}}\right)\cos \varphi_2 }{\sqrt{2(2 + 2\cos \varphi_2 \cos \frac{\gamma}{2} ) \sqrt{1 - \sin^2 \varphi_2  \cos^2 \frac{\gamma}{2} } + 4 - 4\sin^2 \varphi_2 \cos^2 \frac{\gamma}{2} + 2\sqrt{2 + 2\cos \gamma }\cos \varphi_2 }}\,.
\end{multline*}
\end{footnotesize}
\ls{Using MAPLE we obtain that the minimum value of $\rho^\b_{\gamma}(\varphi_2)$ is zero with respect to $\gamma \in (0 \,,\, 2\pi/5]$ and $\varphi_2 \in [0 \,,\, 2\pi)$.} This implies the thesis.
\end{proof}

Note that the minimum value of $\rho^\b_{\gamma}(\varphi_2)$ is achieved for $\gamma = 2\pi/5$ and $\varphi_2 = \pi$; it corresponds to the vanishing of $\mathbf{I}_\gamma(\varphi_2)$, see also Proposition~\ref{propB}.
\begin{figure}[t!]
\begin{center}
    \includegraphics[width=0.33\linewidth]{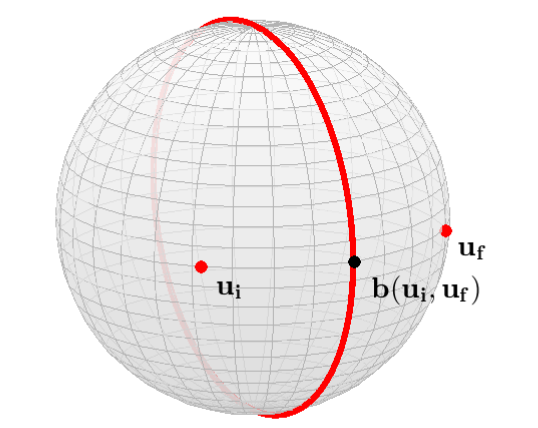}
    \includegraphics[width=0.33\linewidth]{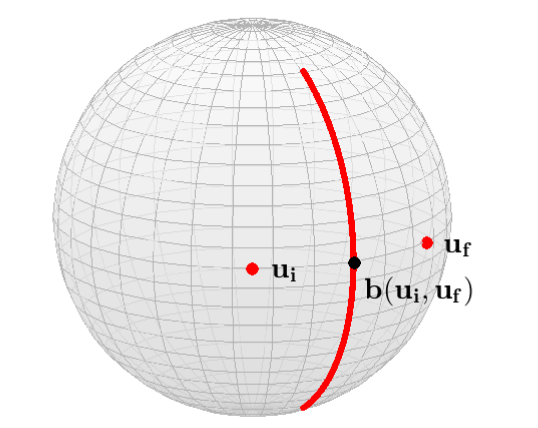}
    \caption{Left: the locus of point described by $\mathbf{S}_{\gamma}(\varphi_2)$ (red line) for $\varphi_2 \in [0 \,,\, 2\pi)$ and $\gamma = \pi/2 > 2\pi/5$, coincident with $\BC(\u_i,\u_f)$. Right: the same locus is described for $\gamma = \pi/3 <2\pi/5$. The initial and final tangent, $\u_i$ and $\u_f$ (red dots), as well as their bisector $\bis$ (black dot), are also shown.}
    \label{fig:S_directions}
\end{center}
\end{figure}
Figure~\ref{fig:S_directions} shows the portion of $\BC(\s_0,\s_4)$ covered by $\S_\gamma(\varphi_2)$ for two different values of $\gamma$. Finally, the next corollary is of specific interest for our local interpolation algorithm,
\begin{corollary} \label{1sufficientcondition}
Let the unit displacement $\mathbf{\Delta} \u$ in $\mathbb{S}_2$ be coincident with any spherical point belonging to $\BC(\s_0,\s_4).$ Then,  if $\gamma > 2\pi/5,$ there exists $ \varphi_2 \in [0\,,\,2\pi)$ such that the unit PH displacement $\S_\gamma(\varphi_2)$ coincides with $\mathbf{\Delta} \u.$
\end{corollary}
\begin{proof}
Under the considered hypothesis, Proposition \ref{propF} states that $\S_\gamma(\varphi_2)$ covers the full circle $\BC(\s_0,\s_4)$ when $\varphi_2$ varies in $[0\,,\,2\pi).$ This directly imples the thesis. \end{proof}

\begin{remark}\label{empirical_considerations}
Note that in all our experiments, when $\gamma > 2\pi/5,$ then $\S_\gamma(\varphi_2)$ continuously \ls{covers} $\BC(\s_0,\s_4)$ exactly once. Conversely, when $\gamma < 2\pi/5,$  for $\varphi_2 \in [0,2\pi],$ starting form $\b$ at $\varphi_2 = 0, \S_\gamma(\varphi_2) $ \ls{covers} a portion of $\BC_h(\s_0,\s_4)$ exactly \ls{twice}. In the special case $\gamma = 2\pi/5,$ for $\varphi_2 \in [0,\pi) ((\pi,2\pi]),$ it \ls{covers} exactly once the quarter of $\BC(\s_0,\s_4)$ in the quadrant generated by $\b$ and $\cross$ ($-\cross$), with 
$$ \lim_{\varphi_2 \rightarrow  \pi ^\mp} \S_{2\pi/5}(\varphi_2) =  \pm \ \cross.$$ 
%Proving this analytically was not possible however, even relying on the help of \texttt{MAPLE}.     
\end{remark}

%Whenever it is possible to find $\varphi_2$ s.t. $\S_\gamma( \varphi_2) = \mathbf{\Delta}_\u$ (for example when  %the hypothesis of Corollary~\ref{1sufficientcondition} holds true), then the free positive parameter $\mu$ has %to be chosen as in (\ref{lamval}) in order to fully satisfy \eqref{alg_integral}. 

\begin{figure}[t!]    
\begin{center}

    \includegraphics[width=0.49\linewidth]{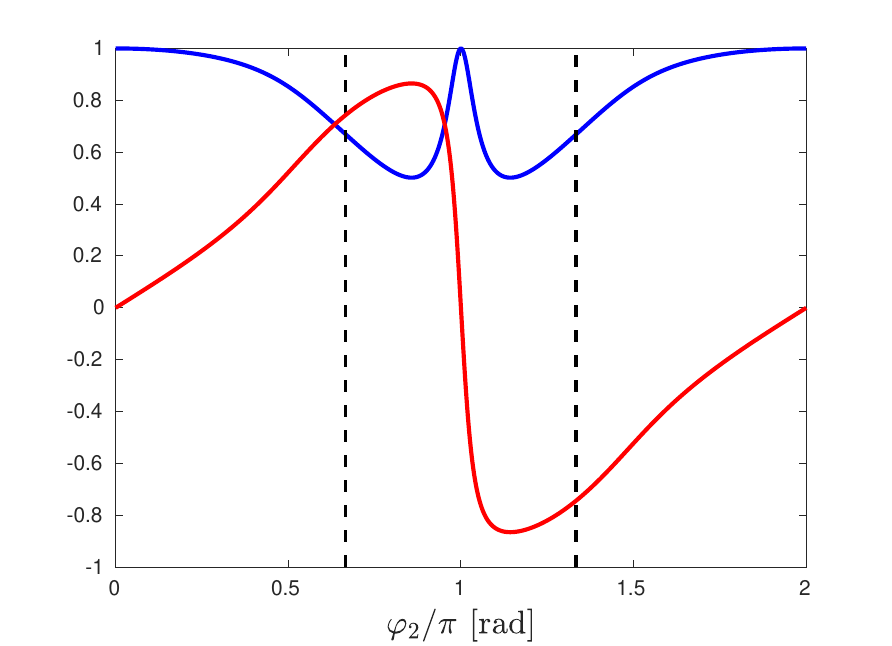}
    \includegraphics[width=0.49\linewidth]{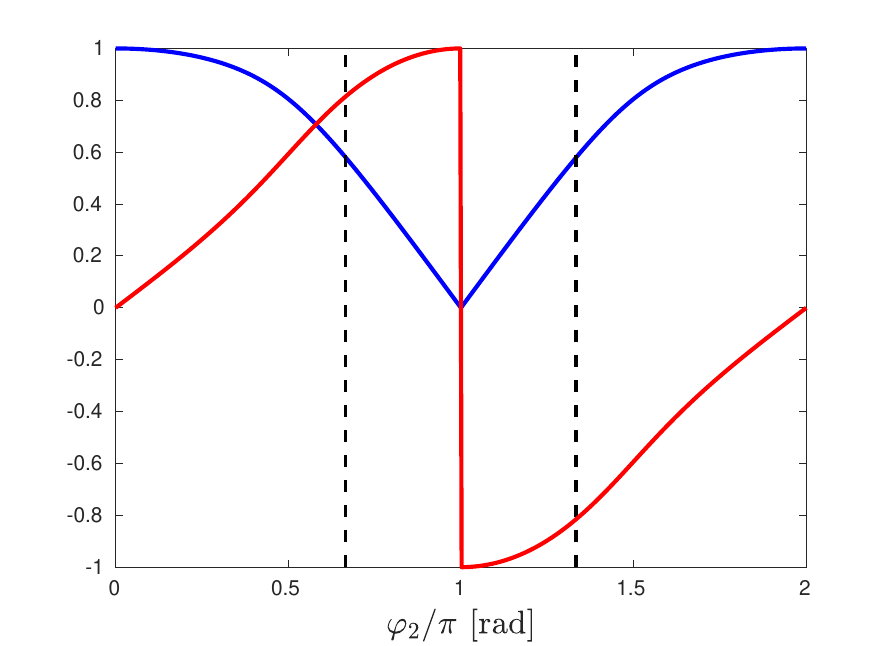}
    \includegraphics[width=0.49\linewidth]{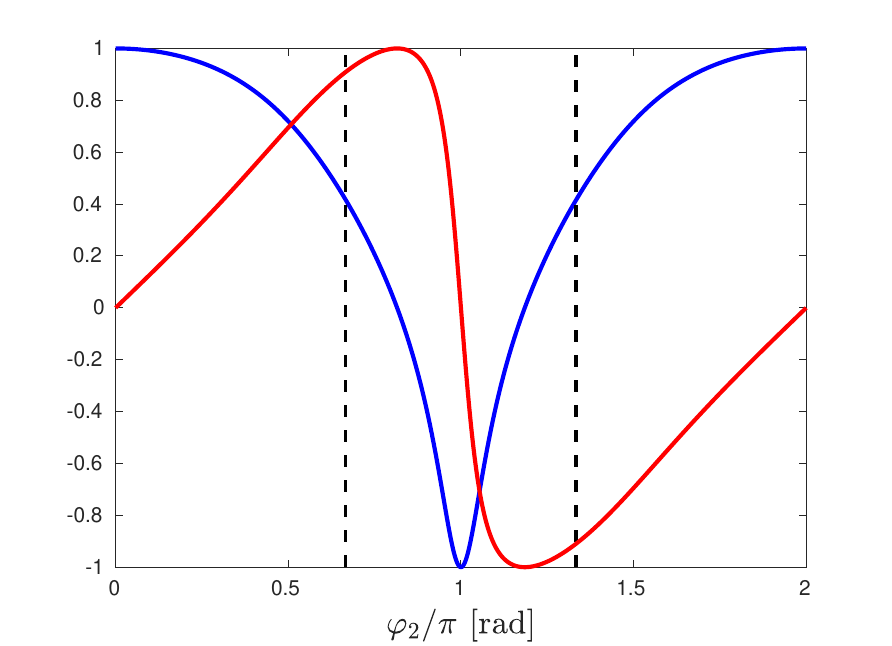}
    \caption{The blue lines represent $\S_{\gamma}(\varphi_2) \cdot \bis$, the red lines $\S_{\gamma}(\varphi_2) \cdot \cross$, for $\gamma = \pi/3$ (top left), $\gamma = 2\pi/5$ (top right), $\gamma = \pi/2$ (bottom). The value \ls{$ \varphi_2/\pi = 2/3$} (and the symmetric one \ls{$\varphi_2/\pi = 4/3$}) relative to the sufficient condition~\eqref{suffcond} is depicted with a black dashed line.}
    \label{bisnorm}  
\end{center}
\end{figure}

Figure~\ref{bisnorm} shows an example of these three cases. We conclude with a sufficient condition that will be used in the spline extension of the method presented in Section~\ref{sec:global_alg} to ensure the presence of solutions even if $\gamma < 2\pi/5$.
\begin{proposition} \label{prop:suffcond}
A sufficient condition for the existence of at least \ls{one} RRMF5-I curve $\r$ solution  to the interpolation problem described by equations \eqref{C0interp}--\eqref{G1interp} is that
\begin{equation}\label{suffcond}
    \b \cdot \mathbf{\Delta} \u > \b \cdot \mathbf{S}_{\gamma} \left(\frac{2 \pi}{3}\right).   
\end{equation}
\end{proposition}
\begin{proof}
As stated in Proposition~\ref{propB}, it holds that
\[
   \mathbf{S}_{\gamma} (0) =  \b.
\]
If condition \eqref{suffcond} is met, since $\mathbf{I}_{\gamma}(\varphi_2)$ does not vanish in the interval $[0\,,\, 2 \pi/3]$, for continuity it will exist a parameter value $\varphi_2 $ so that
\[
    \mathbf{\Delta} \u = \mathbf{S}_{\gamma} (\varphi_2).
\]
\end{proof}

\begin{remark} Instead of using the reference value $\varphi_2 = 2 \pi/3$ considered in Proposition \ref{prop:suffcond}, it could be preferable to determine the analytic expression of the minimum of the  function $\b \cdot \mathbf{S}_{\gamma}(\varphi_2).$  Even if this was not possible because of the complexity of $\mathbf{S}_{\gamma}(\varphi_2)$ in terms of $\gamma$ and $\varphi_2$,  we empirically observed  that, for  different values of $\gamma \in (0\,,\,\pi),$ the value $\varphi_2 = 2\pi/3$ is usually close to the minimum. 
\end{remark}

\subsection{The bisection algorithm}\label{sec:bisection}
When the existence of one value (or possibly two values)  $\hat \varphi_2 = \hat \varphi_2 (\gamma)$ such that  $\mathbf{S}_{\gamma}(\hat \varphi_2) = \mathbf{\Delta} \u$ is ensured by Proposition~\ref{propF} or Proposition~\ref{prop:suffcond}, in view of the complexity of the derivation for the analytic expression of $\hat \varphi_2$, we rely on a bisection algorithm to find a suitable approximation of its value. In order to explain the related details, let us first note that 
\begin{itemize}
    \item $\mathbf{\Delta}\u \cdot \cross > 0$ requires $\hat \varphi_2 \in (0 \,,\, \pi)$;
    \item $\mathbf{\Delta}\u \cdot \cross < 0$ requires $\hat \varphi_2 \in (-\pi \,,\, 0)$;
    \item $\mathbf{\Delta}\u \cdot \cross = 0$ requires $\hat \varphi_2 = 0$ or $\hat \varphi_2 = \pi$.
\end{itemize}  
These observations follow from Proposition~\ref{propN}, see also Figure~\ref{bisnorm}.
Let us then start by checking if the values $\varphi_2 = 0$ or $\varphi_2 = \pi$ ensure $\S_{\gamma}(\varphi_2) = \mathbf{\Delta}\u.$ If this is not the case, we choose between the positive or negative range for $\varphi_2$ considering the sign of $\mathbf{\Delta}\u \cdot \cross.$ Let us refer for example to the case $(0\,,\,\pi).$ Then, assuming that either the hypothesis of Proposition~\ref{propF} or the sufficient condition \eqref{suffcond} hold\ls{s} true, the bisection algorithm to find a zero of 
\begin{equation} \label{bisection_fun}
    f(\varphi_2) := \left(  \S_{\gamma}(\varphi_2) - \mathbf{\Delta}\u \right) \cdot \b  
\end{equation}
is suitably exploited.
From Proposition~\ref{propB}, $f(0) > 0$. We thus have three possibilities: 
\begin{itemize}
\item if $\gamma > 2 \pi/5$, from Proposition~\ref{propB} it holds that $f(\pi) < 0$. Thus $I = [0\,,\, \pi]$ can be considered as initial interval for a bisection algorithm to approximate \ls{a} solution of $f(\varphi_2) = 0$ in this range.
\item if $\gamma = 2 \pi/5$ we assume the fulfillment of the sufficient condition~\eqref{suffcond}. Thus, $f(2\pi/3) < 0$ which implies that $I = [0\,,\, 2\pi/3]$ can be considered as initial interval for a bisection algorithm to approximate \ls{a} solution in $(0\,,\,\pi)$ of $f(\varphi_2) = 0$.
\item if $\gamma < 2 \pi/5$, from Proposition~\ref{propB} we know that $f(\pi) > 0.$  On the other hand, since we are assuming that the sufficient condition~\eqref{suffcond} is fulfilled, we also have  $f(2\pi/3) < 0$. Thus, in this case, both $I_1 = [0\,,\, 2\pi/3]$ and $I_2 = [2\pi/3\,,\,\pi]$ can be considered as initial intervals for a bisection algorithm to approximate two \ls{distinct} solutions of $f(\varphi_2) = 0$ in  $(0\,,\,\pi)$. Calling $S_{\varphi_2}$ the set of the two resulting values, we choose the best one solving
\begin{equation*} 
\argmin_{\varphi_2 \in S_{\varphi_2}} \left( g(\varphi_2) \right), \qquad \text{where}  \qquad g(\varphi_2) := \sum_{i=0}^3 \arccos(\s_i \cdot \s_{i+1})
\end{equation*}
is a functional measuring the angular amplitude of the control polygon of the tangent indicatrix.
\end{itemize}
We observe that the empirical analysis reported in Remark~\ref{empirical_considerations} suggests that, for $\gamma > 2\pi/5$, the solution of $f(\varphi_2) = 0$ belonging to $(0,\pi)$ is unique. It also suggests that exactly two (one) solutions of $f(\varphi_2) = 0$ exist in $(0,\pi)$, when \eqref{suffcond} is fulfilled and $\gamma < 2 \pi /5 \,(\gamma = 2 \pi /5)$.

\section{Global spline extension} \label{sec:global_alg}

We now present the spline extension of the local $G^1$ RRMF5-I interpolation algorithm previously introduced. The global algorithm takes as input data a stream of points $\p_k \in \mathbb{E}^3$, for $k = 0,\ldots,N$, with $\p_k \neq \p_{k+1}$, together with an initial frame orientation $(\u_0,\mathbf{v}_0,\w_0).$  Note that the spline segments are locally constructed one after the other in order to obtain from the previous segment the initial frame orientation. 
\subsection{Spline construction}
Let consider
$\p_i^{(k)} := \p_k\,,\,\p_f^{(k)} := \p_{k+1}\,,\,\u_i^{(k)} := \u_k\,,\,\u_f^{(k)} :=\u_{k+1}, \, \Delta \p^{(k)} := \p_f^{(k)} - \p_i^{(k)}$. In addition, the triples $\left(\u_i^{(k)}, \mathbf{v}_i^{(k)}, \w_i^{(k)}\right)$ and $\left(\u_f^{(k)}, \mathbf{v}_f^{(k)}, \w_f^{(k)}\right)$ refer to the initial and final rotation minimizing frame orientation of the $k$--th spline segment. To ensure the (global) continuity of the frame,  we set 
\[
\left(\u_i^{(k)}, \mathbf{v}_i^{(k)}, \w_i^{(k)} \right) = \left(\u_f^{(k-1)}, \mathbf{v}_f^{(k-1)}, \w_f^{(k-1)}\right)\,,
\]
where ($\u_i^{(0)}$, $\mathbf{v}_i^{(0)}$, $\w_i^{(0)}$) is chosen as the assigned triple $(\u_0,\mathbf{v}_0,\w_0).$
In view of \eqref{construction_hp}, the right unit tangent $\u_f^{(k)}$ has to be obtained by rotating $\u_i^{(k)}$ about $\mathbf{\Delta}\u^{(k)}\,.$ Setting  $\tau^{(k)} \in (0\,,\,\pi)$ as
    \begin{equation} \label{eq:tau}
    \tau^{(k)} := \arccos(\u_i^{(k)} \cdot \mathbf{\Delta}\u^{(k)})\,,
\end{equation}
with $\mathbf{\Delta}\u^{(k)} :=  \Delta \p^{(k)} /\,\vert\mathbf{\Delta}\p^{(k)} \vert\,,$ we necessarily obtain
\begin{equation} \label{gammamax}
 \gamma^{(k)} := \arccos(\u_i^{(k)}\cdot \u_f^{(k)}) \le \gamma_{max}^{(k)}:= \left\{\begin{array}{ll} 2\tau^{(k)} & \mbox{ if } \tau^{(k)} \le \pi/2\,, \cr  2(\pi-\tau^{(k)}) & \mbox{ if } \tau^{(k)} > \pi/2. \end{array} \right.  
 \end{equation}

The following proposition gives a condition to establish if the interpolation problem considered in the previous section can be solved. 
 \begin{proposition}\label{prop:tangent}
    Given two distinct point $\p_i^{(k)},\, \p_f^{(k)}$, and the triple $\left(\u_i^{(k)}, \mathbf{v}_i^{(k)}, \w_i^{(k)}\right)\,,$ it is possible to define $\u_f^{(k)}$ so that \eqref{construction_hp} holds and an RRMF5-I curve solving \eqref{C0interp}-\eqref{G1interp} exists for these input data if and only if $\tau^{(k)} < 4\pi/5$ .
\end{proposition}
\begin{proof}
 We first prove that, if $\tau^{(k)} \geq 4\pi/5$, it is not possible to define $\u_f^{(k)}$ so that \eqref{construction_hp} holds. From \eqref{gammamax} if follows that $\gamma^{(k)} \leq 2\pi/5$ when $\tau^{(k)} \geq 4\pi/5$ and, consequently, the hypothesis of Proposition~\ref{propF} is not satisfied . Denoting with $\phi^{(k)}$ the angle between $\mathbf{\Delta}\u^{(k)}$ and $\bis(\u_i^{(k)}, \u_f^{(k)}),$ considering some easy angular relations, we can derive that
\begin{equation}\label{eq:tau_gamma_function}
\cos \tau^{(k)} = \cos \frac{\gamma^{(k)}}{2} \ \cos \phi^{(k)}\,.
\end{equation}
Together with \eqref{gammamax}, this implies that $\phi^{(k)} > \pi/2$ when $\tau^{(k)} > \pi/2$ and, as a consequence, $\mathbf{\Delta}\u^{(k)} \not \in \BC_h(\s_0,\s_4)$. In view of Proposition~\ref{prop:gamma}, we can conclude that, when $\tau^{(k)} \ge 4\pi/5$, no selection of the right unit tangent $\u_f^{(k)}$ is admissible for the local algorithm.   
Note however that this situation corresponds to require almost a full reversion of the motion direction, see Example~\ref{ex:badstream} below.

For what concern the second part of the proof, we preliminary note that, if  $\u_f^{(k)}$ is such that $ \gamma^{(k)} = \gamma_{max}^{(k)}\,,$ then the three spherical points $\u_i^{(k)}, \u_f^{(k)}$ and $\mathbf{\Delta}\u^{(k)}$ are  coplanar, with $\mathbf{\Delta}\u^{(k)}$ coincident \ls{with $\bis(\u_i^{(k)}, \u_f^{(k)})$ if $\tau^{(k)} <  \pi/2$ and with $-\bis(\u_i^{(k)}, \u_f^{(k)})$ if $\tau^{(k)} >  \pi/2$.} As a consequence, if $\tau^{(k)}< \pi /2$ there exists a choice of $\u_f^{(k)}$ satisfying condition~\eqref{suffcond}. When  $\pi/2 \leq \tau^{(k)} < 4 \pi/5\,,$ \eqref{gammamax} implies that there is a choice of  $\u_f^{(k)}$ such that $ \gamma^{(k)} > 2 \pi /5\,.$ Thus Proposition~\ref{propF} implies the thesis.  
\end{proof}

Let us now introduce our strategy for defining the right unit tangent, limiting ourselves to the case $\tau^{(k)} < 4\pi/5.$ We start with a preliminary computation of a reference unit tangent $\u^{\text{ref}}_{k+1}$, not necessarily fulfilling the condition in \eqref{construction_hp} with $\u_i^{(k)}$ and $\Delta \p^{(k)}$. The missing right unit tangent $\u_f^{(k)}$ is derived solving the local optimization problem, 
\begin{equation} \label{tangentgen}
\u_f^{(k)} = \max_{\u \in S^{(k)}} ( \u \cdot \u^{\text{ref}}_{k+1})
\end{equation}
where
\begin{footnotesize}
$$ S^{(k)} := \left\{\u \in \mathbb{S}^2: (\u - \u_i^{(k)} )\cdot \mathbf{\Delta}\u^{(k)}=0 \land \left( \bis(\u_i^{(k)},\u) \cdot \left(\mathbf{\Delta}\u^{(k)}-\mathbf{S}^{(k)}_\gamma \left(\frac{2 \pi}{3}\right)\,\right) > 0 \lor \arccos(\u_i^{(k)} \cdot \u) >  \frac{2\pi}{5} \right)\right\}.$$
\end{footnotesize}
We can interpret the solution of the above optimization problem as follows. Among all the possible unit vectors generated by rotating the initial unit tangent $\u_i^{(k)}$ about the unit displacement $\mathbf{\Delta} \u^{(k)},$ we aim at choosing the closest one to the reference unit tangent $\u^{\text{ref}}_{k+1}$ on the sphere, also ensuring the possibility to solve the corresponding local interpolation problem introduced in the previous section. Indeed we require that either the sufficient condition in \eqref{suffcond} is satisfied or the hypothesis of Proposition~\ref{propF} holds true.

If for all the spline segments we can obtain an admissible right unit tangent, it is possible to define a spline path $\X = \X(u)$, with $u \in [0\,,\, u_N]$, so that
$$
\X(u) = \r^{(k)}\left(\frac{u-u_k}{u_{k+1}-u_k}\right)\,, \qquad \mbox{for } u \in [u_k , u_{k+1}], \,\qquad k = 0,\ldots,N-1\,. 
$$
\ls{Here} $\r^{(k)}$ is the $k$--th RRMF5-I curve generated with the algorithm of the previous section, \ls{$u$ is the global parameter and $u_0,\ldots,u_N$ are the spline breakpoints.} 

\ls{It may happen that, if two consecutive displacements, $\bm \Delta \p^{(k-1)}$ and $\bm \Delta \p^{(k)}$ have almost opposite directions, the input data of the $k$--th local interpolation problem do not satisfy the hypothesis of Proposition~\ref{prop:tangent}. As a result, it becomes impossible to compute an admissible right unit tangent. To overcome this problem, it is possible to insert a suitable point $\p_m^{(k)}$ between $\p_f^{(k-1)}\equiv\p_i^{(k)}$ and $\p_f^{(k)}\equiv\p_i^{(k+1)}$ to reduce the angular distance between successive displacements. Algorithm~\ref{alg:point_insertion} presents a simple procedure which requires in input a free parameter $c\in(0,1]$ and returns $\p_m^{(k)}$ and a corresponding unit tangent $\u_m^{(k)}$. The parameter $c$ is simply used in step 1 of the algorithm to identify an intermediate necessary point $\p_c^{(k)}$ as convex combination of $\p_i^{(k)}$ and $\p_f^{(k)}$, see also Example~\ref{ex:badstream}. The possibility of inserting intermediate Hermite data according to the procedure described in Algorithm~\ref{alg:point_insertion} guarantees the existence of RRMF5-I spline interpolants of arbitrary 3D data stream.}

\bigskip

\ls{\begin{algorithm}[H]
\KwIn{The points $\p_i^{(k)}$, $\p_f^{(k)}$, the unit tangent $\u_i^{(k)}$, a parameter $c$, with $0 < c \leq 1$.}
\KwOut{The point $\p_m^{(k)}$, the unit tangent $\u_m^{(k)}$.}
Compute the point $\p_c^{(k)} := (1-c)\,\p_i^{(k)}+c\,\p_f^{(k)}$\;
Compute the unit displacement $\mathbf{\Delta} \u^{(k)}=\mathbf{\Delta} \p^{(k)}/\vert\mathbf{\Delta} \p^{(k)}\vert\,,$ with $\mathbf{\Delta} \p^{(k)}=\p_f^{(k)}-\p_i^{(k)}$\;
Compute the unit vector $\u_c^{(k)} := (\mathbf{\Delta} \u^{(k)}\times\u_i^{(k)})\times \mathbf{\Delta} \u^{(k)}/\vert\mathbf{\Delta} \u^{(k)}\times\u_i^{(k)} \vert$\;
Compute the unit bisector $\bis(\u_i^{(k)},\mathbf{\Delta} \u^{(k)})$ according to \eqref{bdef}\;
Find $(r^*\,,\,s^*) \in \mathbb{R}^2$ such that $\p_i^{(k)} + r \,\bis(\u_i^{(k)},\mathbf{\Delta} \u^{(k)})= \p_c^{(k)}+s\,\u_c^{(k)}$\;
Compute $\p_m^{(k)} := \p_c^{(k)}+s^*\,\u_c^{(k)}$\;
Set $\u_m^{(k)} := \mathbf{\Delta} \u^{(k)}$\;
\caption{Insertion of the point $\p_m^{(k)}$}
\label{alg:point_insertion}
\end{algorithm}}

\bigskip
Algorithm~\ref{alg:point_insertion}  identifies the new point $\p_m^{(k)}$ to be inserted as intersection between the coplanar straight lines in the direction of the bisector $\bis(\u_i^{(k)},\mathbf{\Delta} \u^{(k)})$ and of the unit vector $\u_c^{(k)}$ defined in step~3. Then two new local interpolation subproblems replace the original $k$--th problem. The first one, identified here below by the superscript $k_0$ has the Hermite data $\p_i^{(k_0)}= \p_i^{(k)},\,\p_f^{(k_0)}= \p_m^{(k)},\,\u_i^{(k_0)}=\u_i^{(k)},\,\u_f^{(k_0)}=\u_m^{(k)}$, as shown on the left of Figure~\ref{fig:addpoint}. We can easily observe that $ \tau^{(k_0)}= \tau^{(k)}/ \,2$ by construction. Dealing with planar data, it follows from \eqref{eq:tau_gamma_function} that $\gamma^{(k_0)} = 2 \,\tau^{(k_0)}$, and, as a result, $\gamma^{(k_0)} = \tau^{(k)}> 2\pi/5$. The hypothesis of Corollary~\ref{1sufficientcondition} is then satisfied and, as a consequence, the corresponding local problem admits a solution. For the second local interpolation subproblem, identified by the superscript $k_1$, instead we have $\p_i^{(k_1)}= \p_m^{(k)},\,\p_f^{(k_1)}= \p_f^{(k)},\,\u_i^{(k_1)}=\u_m^{(k)}$, as shown on the right of Figure~\ref{fig:addpoint}.  Since $\tau^{(k_1)}<\pi/2$, the unit vector $\u_f^{(k_1)}$ can be derived from formula given in \eqref{tangentgen}. 

\begin{figure}[t!]
\centering\begin{footnotesize}
\begin{subfigure}[b]{0.5\textwidth}
\centering
\begin{tikzpicture}
    \coordinate (pi) at (0,0);
    \coordinate (pf) at (5,0);
    \coordinate (A) at (1,0);
    \coordinate (pm) at (1,3);
    \coordinate (B) at (1.3,3.9);
    \coordinate (um) at (3,3);
    \coordinate (u_i) at (-1.5,1.125);
    \draw[thick] (pi) -- (pf);
    \draw[draw=blue,->, thick] (pi) -- (u_i); 
    \draw[draw=blue,thick] (pi) -- (pm); 
    \draw[draw=blue,->,thick] (pm) -- (um); 
    \draw[draw=black, dashed,thin] (pm) -- (B); 
    \draw[draw=black, dashed, thin] (A) -- ++(0,4);  
    \pic [draw=black, angle radius=0.4cm,thin] {angle = pm--pi--u_i};
    \pic [draw=black, angle radius=0.4cm,thin] {angle = pf--pi--pm};
    \pic [draw=black, angle radius=0.4cm,thin] {angle = um--pm--B};    
    \node at ($(pi)+(1,0)$) [below left,text=blue] {$\bm{p}_i^{(k_0)}$ \textcolor{black}{= $\bm{p}_i^{(k)}$}};
    \node at (pf) [below right] {$\bm{p}_f^{(k)}$};
    \node at (pm) [above left,text=blue] {$\bm{p}_f^{(k_0)}$ \textcolor{black}{= $\bm{p}_m^{(k)}$}};
    \node at ($(pi)!0.5!(u_i)-(0.8,0)$)  [below,text=blue] {$\bm{u}_i^{(k_0)}$ \textcolor{black}{= $\bm{u}_i^{(k)}$}}; 
    \node at ($(pi)!0.5!(pf)$)  [below] {$\bm{\Delta} \bm{p}^{(k)}$}; 
    \node at ($(pi)!0.5!(pm)+(0.1,0)$)  [above left,text=blue] {$\bm{\Delta} \bm{p}^{(k_0)}$}; 
    \node at ($(pm)!0.5!(um)$)  [below,text=blue] {$\bm{u}_f^{(k_0)}$ \textcolor{black}{= $\bm{u}_m^{(k)}$}}; 
    \node at (-0.65,0.45) [above right] {$\tau^{(k_0)}$}; 
    \node at (1.35,3.25) [above right] {$\tau^{(k_0)}$}; 

    \filldraw[red] (pi) circle (2pt); 
    \filldraw[red] (pf) circle (2pt); 
    \filldraw[blue] (pm) circle (2pt);
\end{tikzpicture}
\end{subfigure}
\hspace*{-0.5cm}
\begin{subfigure}[b]{0.5\textwidth}
\centering
\begin{tikzpicture}
    \coordinate (pi) at (0,0);
    \coordinate (pf) at (5,0);
    \coordinate (A) at (1,0);
    \coordinate (pm) at (1,3);
    \coordinate (B) at (1.3,3.9);
    \coordinate (um) at (3,3);
    \coordinate (u_i) at (-1.5,1.125);

    \draw[thick] (pi) -- (pf);
    \draw[->, thick] (pi) -- (u_i); 
    \draw[dashed,thin] (pi) -- (pm); 
    \draw[draw=blue,->,thick] (pm) -- (um); 
    \draw[draw=blue,thick] (pm) -- (pf); 
    \draw[draw=black, dashed,thin] (pm) -- (B); 
    \draw[draw=black, dashed, thin] (A) -- ++(0,4);  

    \pic [draw=black, angle radius=0.6cm,thin] {angle = pf--pm--um};
    \pic [draw=black, angle radius=0.6cm,thin] {angle = pm--pf--pi};
  
    \node at (pi) [below left] {$\bm{p}_i^{(k)}$};
    \node at ($(pf)-(1,0)$) [below right,text=blue] {$\bm{p}_f^{(k_1)}$  \textcolor{black}{= $\bm{p}_f^{(k)}$}};
    \node at (pm) [below left,text=blue] {$\bm{p}_i^{(k_1)}$ \textcolor{black}{= $\bm{p}_m^{(k)}$}};
    \node at ($(pi)!0.5!(u_i)-(0.2,0)$)  [below] {$\bm{u}_i^{(k)}$}; 
    \node at ($(pi)!0.5!(pf)$)  [below] {$\bm{\Delta} \bm{p}^{(k)}$}; 
    \node at ($(pm)!0.5!(pf)$)  [above right,text=blue] {$\bm{\Delta} \bm{p}^{(k_1)}$}; 
    \node at ($(pm)!0.5!(um)+(0.1,0)$)  [above,text=blue] {$\bm{u}_i^{(k_1)}$ \textcolor{black}{= $\bm{u}_m^{(k)}$}}; 
    \node at ($(pm)+(0.7,-0.55)$) [above right] {$\tau^{(k_1)}$}; 
    \node at ($(pf)+(-1.55,+0.05)$) [above right] {$\tau^{(k_1)}$}; 
    
    \filldraw[red] (pi) circle (2pt); 
    \filldraw[red] (pf) circle (2pt); 
    \filldraw[blue] (pm) circle (2pt);
\end{tikzpicture}
\end{subfigure}\end{footnotesize}
\caption{\ls{Insertion of the additional point $\p_ m^{(k)}$ according to Algorithm~\ref{alg:point_insertion}. The two new local subproblems are shown with index $k_0$ (left) and $k_1$ (right).}}
\label{fig:addpoint}
\end{figure}

\subsection{Numerical examples}
We now present a selection of numerical experiments to highlight the performance of our approach. We use point data streams, either (i) sampled from analytical spatial curves or (ii) freely assigned to describe generic 3D paths. 

For each example of type (i) we consider an analytical curve $\C(u)$, with $u \in [0,U]$, and we set a uniform parametric grid
$$
u_k = k\, \frac{U}{N}, \quad \text{with} \quad k = 0,\ldots,N\,,
$$
to sample points and reference unit tangents
$$
\p_k = \C\left(u_k\right), \qquad 
\u^{\text{ref}}_k = \frac{\C'\left(u_k\right)}{\vert\C'\left(u_k\right)\vert}.
$$
\begin{figure}[t!]
\begin{center}
    \includegraphics[width=0.32\linewidth]{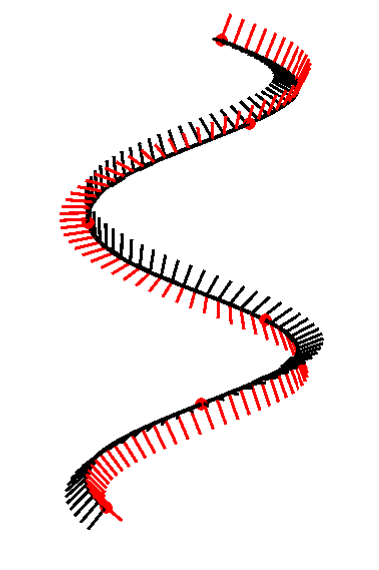}
    \includegraphics[width=0.32\linewidth]{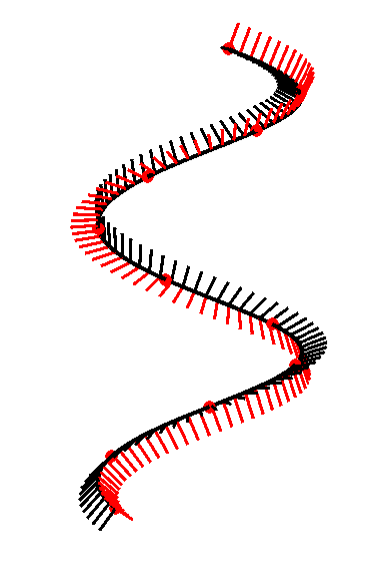}
    \includegraphics[width=0.32\linewidth]{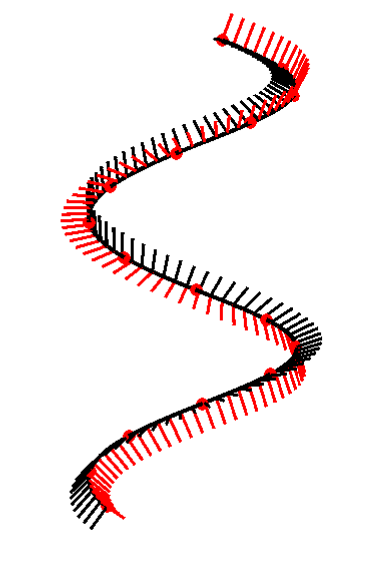}
    \caption{The analytic curve of Example~\ref{ex:helix} is depicted (black dashed lines), together with the RRMF5-I splines (black lines) and their attached frames (black and red arrows), resulting from the interpolation of 6 (left), 11 (center), 16 (right) sampled data points (red dots).}  
    \label{fig:helix}
\end{center}
\end{figure}
\begin{example}\label{ex:helix} \mbox{}\rm
    Let us consider three distinct sequences of 3D data points sampled from a circular arc-length parameterized helix
\begin{equation*}
\C(u) = \left(\begin{array}{l} 
10 \sin\left( \frac{u}{u_h} \right) \cr
10 \cos\left( \frac{u}{u_h} \right) \cr
-4 \frac{u}{u_h}\cr
\end{array} \right) \qquad \mbox{with } u_h = 2\sqrt{29}, \qquad u \in[0 \,,\,  3.6 \pi u_h]\,.
\end{equation*}  
The resulting RRMF5-I splines are shown in Figure~\ref{fig:helix}.
\end{example}

\begin{figure}[t!]
\begin{center}
    \includegraphics[width=0.49\linewidth]{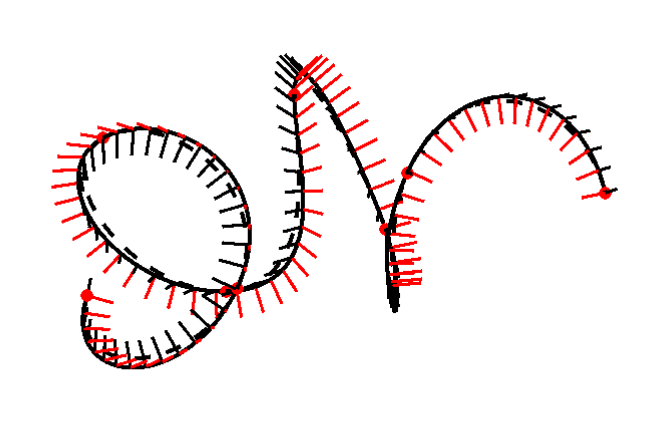}
    \includegraphics[width=0.49\linewidth]{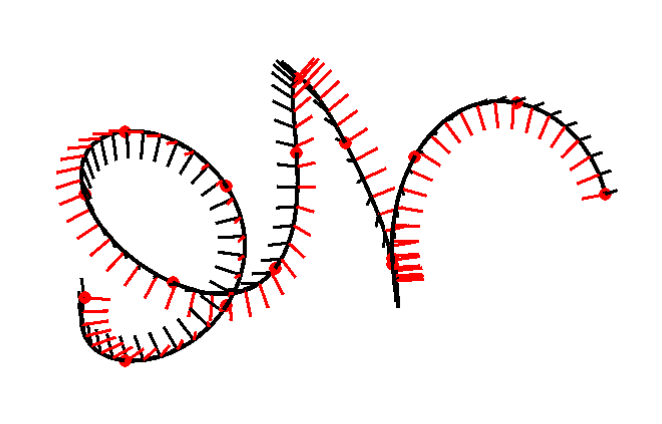}
    \caption{The analytic curve of Example~\ref{ex:torus} is depicted (black dashed lines), together with the RRMF5-I splines (black lines) and their RMF vectors (black and red arrows), resulting from the interpolation of 8 (left) and 16 (right) sampled data points (red dots).}  
    \label{fig:torus}
\end{center}
\end{figure}

\begin{example}\label{ex:torus} \mbox{}\rm
Let us consider 3D data points which are sampled from a curve on a torus with equations 
\begin{equation*}
\C(u) = \left(\begin{array}{l} 
(20 + 10 \cos(3 u))\cos(\frac{u}{2}) \cr
(20 + 10 \cos(3 u))\sin(\frac{u}{2}) \cr
10 \sin(3 u)\cr
\end{array} \right), \qquad u \in[0 \,,\, 2 \pi]\,.
\end{equation*}
The resulting RRMF5-I splines are shown in Figure~\ref{fig:torus} for two distinct number of samples.
\end{example}

\begin{figure}[t!]
\begin{center}
    \includegraphics[width=0.49\linewidth]{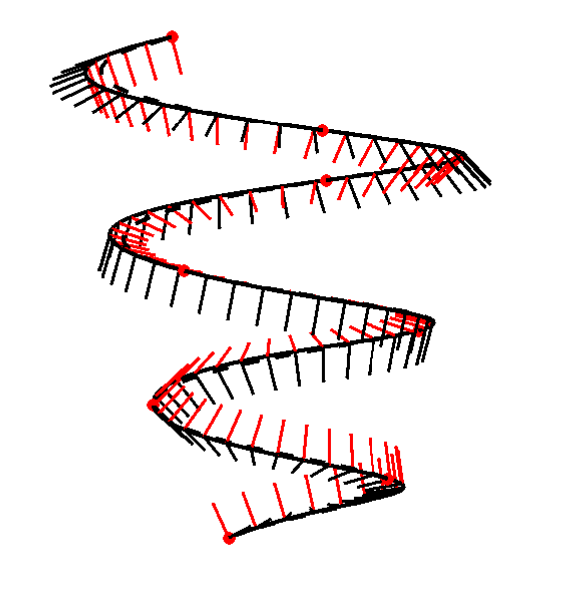}
    \includegraphics[width=0.49\linewidth]{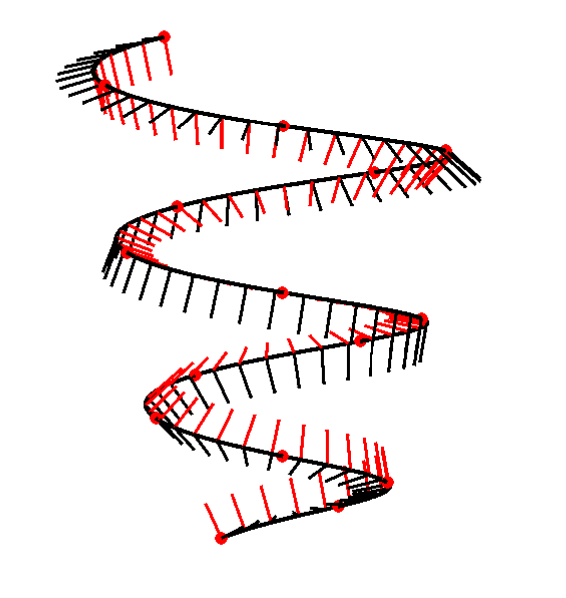}
    \caption{The analytic curve of Example~\ref{ex:spiral} is depicted (black dashed lines), together with the RRMF5-I splines (black lines) and their RMF vectors (black and red arrows), resulting from the interpolation of 8 (left) and 16 (right) sampled data points (red dots).}     
    \label{fig:spiral}
\end{center}
\end{figure}

\begin{example}\label{ex:spiral} \mbox{}\rm
A logarithmic spiral of the form
\begin{equation*}
\C(u) = \left(\begin{array}{l} 
\log(u+3)\, \sin(\pi u) \cr
\log(u+3)\, \cos(\pi u) \cr
\sqrt{u^2 + 4u+ 5} \cr
\end{array} \right), \qquad  u \in\left[0 \,,\,  6\right],
\end{equation*}  
is considered for generating two sequences of 3D data points. The resulting RRMF5-I splines are shown in Figure~\ref{fig:spiral}.
\end{example}

Finally, we present three tests for the application of the $G^1$ RRMF algorithm to generic 3D data stream interpolation. For this second group of examples we set a global chord-length  parameterization as
$$
u_k = u_{k-1} + \Vert \p_k - \p_{k-1}  \Vert
$$ 
with $u_0= 0$. By considering only input point streams, local rules for derivative approximations are used to generate the reference unit tangent $\u^{\text{ref}}_k$.
We rely on the local formulas called MinAJ2 introduced in \cite{Debski1} for data stream applications. More precisely, defining $h_k := u_k - u_{k-1}$, the right derivative $\u^{\text{ref}}_k$ is chosen as
$$
\u_k^{\text{ref}} = \frac{A\,\p_{k-1} + B\,\u_{k-1}^{\text{ref}} + C\,\p_k + D\,\p_{k+1}}{E}, \qquad k = 1,\ldots,N-1,
$$
where
$$
\begin{array}{l} 
A = - h_{k+1}^2(2 h_{k+1}^2 + 6h_{k+1} h_k + 3 h^2_k), \qquad 
B = - h_kh^2_{k+1}(h_{k+1} + h_k)^2, \cr
C = (h_{k+1}+h_k)(2 h_{k+1}^3 + 4 h_{k+1}^2 h_k - h_{k+1} h_k^2 - h_k^3), \cr
D = h^3_k(2h_{k+1} + h_k), 
\qquad E = h_kh_{k+1}(h_{k+1} + h_k)(h_{k+1}^2 + 3h_{k+1} h_k + h^2_k)\,. \cr
\end{array} 
$$
The first and the last references are computed as follows, 
$$
\u^{\text{ref}}_0 = \frac{(\p_1 - \p_0)(h_2+h_1)^2 + (\p_1 - \p_2)h_1^2}{h_1 h_2(h_2 + h_1)}\,,
\qquad
\u^{\text{ref}}_N = -\frac{(\u_{N-1}^{\text{ref}}h_N - 2\p_N + 2\p_{N-1})}{h_N}\,.
$$

\begin{example} \mbox{}\rm \label{ex:badstream}
Consider the sequence of three points 
\[
\p_0 = (0,0,0), \qquad \p_1 = (-5,5,2), \qquad \p_2 = (2,2,0),
\]
for which $\tau^{(1)} = 0.860 \, \pi > 4\pi/5$ and the two displacements have almost opposite direction, i.e., $\arccos(\mathbf{\Delta} \u^{(0)} \cdot \mathbf{\Delta} \u^{(1)}) = 0.863 \, \pi$. We then apply Algorithm~\ref{alg:point_insertion}, considering two different values for the free parameter $c$, namely $c=1/8$ and $c=1/4$. In the first case, we obtain the additional point $\p_m^{(1)} = (-2.2964,8.6050,2.1803)$ to be inserted between $\p_1$ and $\p_2$, together with the unit tangent $\u_m^{(1)} = (0.8890 , -0.3810 , -0.2540)$. In the second case, we have $\p_m^{(1)} = (0.4073,12.2100,2.3605)$, while $\u_m^{(1)}$, which does not depend on $c$, has the same value as before. For both cases, it is then possible to generate the remaining right unit tangent (denoted as $\u_f^{(k_1)}$ in the previous section) leading to the two RRMF5-I spline interpolants shown in Figure~\ref{fig:badstream}. Note that the choice of the parameter $c$ naturally influences the position of $\p_m^{(1)}$. 
\begin{figure}[t!]
\begin{center}
    
    \includegraphics[width=0.45\linewidth]{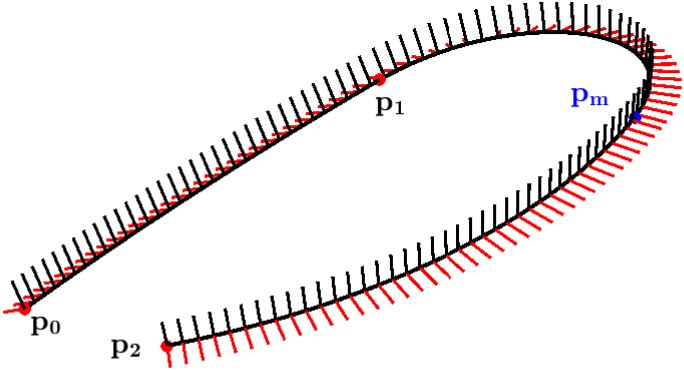}
    \includegraphics[width=0.45\linewidth]{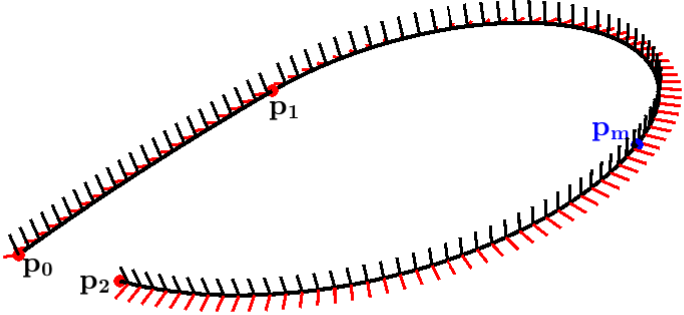}
    \caption{\ls{The RRMF5-I spline curves interpolating the data stream of Example~\ref{ex:badstream} are depicted (black lines), together with their attached frames (black and red arrows), the original interpolation data points $\p_0$, $\p_1$, $\p_2$(red dots) and the additional points $\p_m^{(1)}$ (blue dot) generated by Algorithm~\ref{alg:point_insertion}, with input parameter $c = 1/8$ (left) and $c= 1/4$ (right).}}
    \label{fig:badstream}
\end{center}
\end{figure}
\end{example}

\begin{figure}[t!]
\begin{center}
    \includegraphics[width=0.4\linewidth]{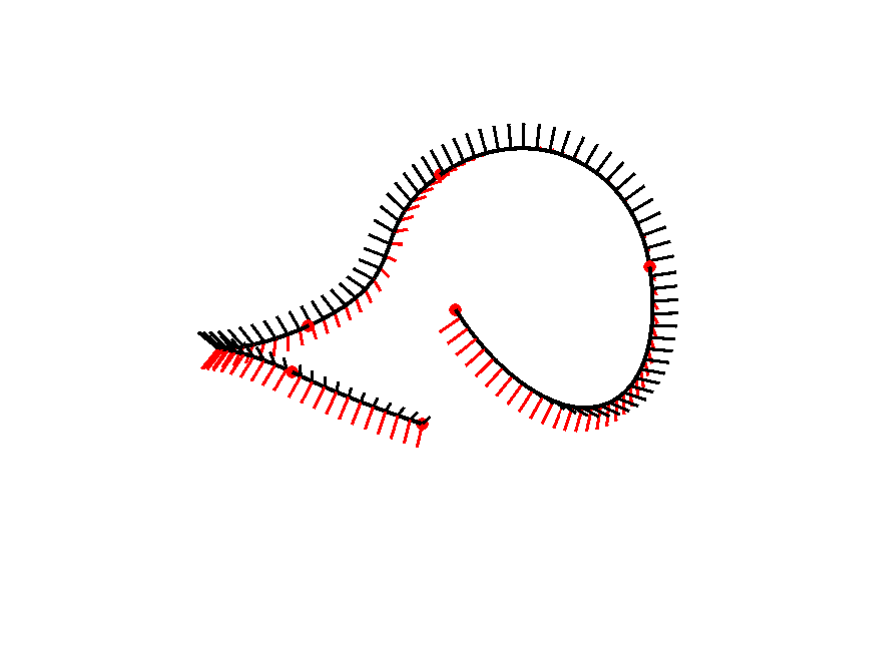}
    \hspace{0.5cm}
    \includegraphics[width=0.4\linewidth]{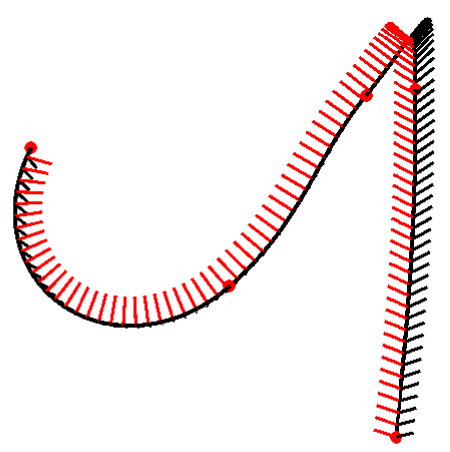}
    \caption{The RRMF5-I spline curves interpolating the data stream of Example~\ref{ex:generic1} (left) and Example~\ref{ex:generic2} (right) are depicted (black lines), together with their attached frames (black and red arrows) and the interpolation data points (red dots).}  
    \label{datastream}
\end{center}
\end{figure}

\begin{example}\label{ex:generic1} \mbox{}\rm
Let us consider the following 3D data stream:
$$
\begin{array}{lll} 
\p_0 = (0 , 0 , 0)^T, & 
\p_1 = (-5 , 5 , 2)^T, & 
\p_2 = (0 , 10 , -2)^T, \cr
\p_3 = (8 , 12 , 5)^T, & \p_4 = (15 , 2 , 3)^T, & \p_5 = (2 , 0 , 7)^T. \cr
\end{array}
$$ 
The resulting RRMF5-I spline is shown in Figure~\ref{datastream} (left).
\end{example}

\begin{example}\label{ex:generic2} \mbox{}\rm
Let us consider the following 3D data stream:
$$
\begin{array}{lllll} 
\p_0 = (0 , 0 , 0)^T, & 
\p_1 = (5 , 5 , 10)^T, & 
\p_2 = (8 , 11 , 9)^T, &
\p_3 = (5 , 14 , 3)^T, & 
\p_4 = (2 , 20 , 7)^T. 
\end{array}
$$ 
The resulting RRMF5-I spline is shown in Figure~\ref{datastream} (right).
\end{example}

\section{Conclusion}\label{sec:conclusions}
An alternative approach to the characterization of RRMF quintic curves of class I, based on the analysis of the projection of the hodograph control points on the unit sphere, was presented. This geometric characterization directly leads to a novel algorithm for the generation of piecewise RRMF motions. In particular, through the request of the symmetric constraint on the input data for the considered local $G^1$ interpolation problem, it was possible to derive sufficient conditions that ensure the existence of solutions. Given a 3D data stream, an algorithm to provide admissible tangents for the construction of a piecewise RRMF spline curve was derived. A simple condition indicates the existence of an admissible sequence of tangents that guarantee the existence of the solution to any local interpolation problem to be considered. It was also shown that, when this is not guaranteed, it is always possible to add an intermediate point \ls{using a simple algorithm which enables the generation of the desired RRMF motion for arbitrary 3D input data stream.} 

\section*{Acknowledgements}
AS and CG acknowledge the contribution of the National Recovery and Resilience Plan, Mission 4 Component 2 – Investment 1.4 -NATIONAL CENTER FOR HPC, BIG DATA AND QUANTUM COMPUTING – funded by the European Union – NextGenerationEU – (CUP B83C22002830001). 
The partial support of the Italian Ministry of University and Research (MUR) through the PRIN projects COSMIC (No. 2022A79M75) and NOTES (No. P2022NC97R), funded by the European Union - Next Generation EU, is also acknowledged. 
CG and AS are members of the INdAM Research group GNCS. The INdAM support through GNCS (CUP E53C24001950001) is gratefully acknowledged.

\bibliographystyle{amsplain}

%%%%%%%%%%%%
\bibliography{biblio}

\providecommand{\bysame}{\leavevmode\hbox to3em{\hrulefill}\thinspace}
\providecommand{\MR}{\relax\ifhmode\unskip\space\fi MR }
\providecommand{\MRhref}[2]{%
  \href{http://www.ams.org/mathscinet-getitem?mr=#1}{#2}
}
\providecommand{\href}[2]{#2}
\begin{thebibliography}{10}

\bibitem{bishop75}
Richard Bishop, \emph{{There is More than One Way to Frame a Curve}}, \AMM
  \textbf{82} (1975), 246--251.

\bibitem{ch02}
Hyeong~In Choi and Chang~Yong Han, \emph{{Euler–Rodrigues} frames on spatial
  {Pythagorean-hodograph} curves}, \CAGD \textbf{19} (2002), 603--620.

\bibitem{choi02}
Hyeong~In Choi, Doo~Seok Lee, and Hwan~Pyo Moon, \emph{{Clifford Algebra, Spin
  Representation, and Rational Parameterization of Curves and Surfaces}}, \ACM
  \textbf{17} (2002), 5--48.

\bibitem{Debski1}
Roman Debski, \emph{Real-time interpolation of streaming data}, \CS \textbf{21}
  (2020), 515--534.

\bibitem{f2010}
Rida~T. Farouki, \emph{{Quaternion and Hopf map characterizations for the
  existence of rational rotation-minimizing frames on quintic space curves}},
  \ACM \textbf{33} (2010), 331--348.

\bibitem{farouki16}
\bysame, \emph{Rational rotation-minimizing frames-recent advances and open
  problems}, \AMC \textbf{272} (2016), 80--91.

\bibitem{farouki02}
Rida~T. Farouki, Mohammad al~Kandari, and Takis Sakkalis, \emph{Structural
  invariance of spatial {Pythagorean hodographs}}, \CAGD \textbf{19} (2002),
  395--407.

\bibitem{fggss17}
Rida~T. Farouki, Graziano Gentili, Carlotta Giannelli, Alessandra Sestini, and
  Caterina Stoppato, \emph{A comprehensive characterization of the set of
  polynomial curves with rational rotation-minimizing frames}, \ACM \textbf{43}
  (2017), 1--24.

\bibitem{fgms09}
Rida~T. Farouki, Carlotta Giannelli, Carla Manni, and Alessandra Sestini,
  \emph{Quintic space curves with rational rotation-minimizing frames}, \CAGD
  \textbf{26} (2009), 580--592.

\bibitem{fgms2012}
\bysame, \emph{{Design of rational rotation–minimizing rigid body motions by
  Hermite interpolation}}, \MC \textbf{81} (2012), 879--903.

\bibitem{review19}
Rida~T. Farouki, Carlotta Giannelli, and Alessandra Sestini, \emph{New
  developments in theory, algorithms, and applications for
  {Pythagorean}--hodograph curves}, Advanced Methods for Geometric Modeling and
  Numerical Simulation (Carlotta Giannelli and Hendrik Speleers, eds.),
  Springer International Publishing, 2019, pp.~127--177.

\bibitem{fhds13}
Rida~T. Farouki, Chang~Yong Han, Petroula Dospra, and Takis Sakkalis,
  \emph{{Rotation-minimizing Euler-Rodrigues rigid-body motion interpolants}},
  \CAGD \textbf{30} (2013), 653--671.

\bibitem{fs12}
Rida~T. Farouki and Takis Sakkalis, \emph{A complete classification of quintic
  space curves with rational rotation-minimizing frames}, \JSC \textbf{47}
  (2012), 214--226.

\bibitem{gss22}
Carlotta Giannelli, Lorenzo Sacco, and Alessandra Sestini, \emph{{
  Interpolation of 3D data streams with $C^2$ PH quintic splines}}, \ACM
  \textbf{48} (2022), 61.

\bibitem{guggenheimer89}
Heinrich Guggenheimer, \emph{Computing frames along a trajectory}, \CAGD
  \textbf{6} (1989), 77--78.

\bibitem{han08}
Chang~Yong Han, \emph{Nonexistence of rational rotation-minimizing frames on
  cubic curves}, \CAGD \textbf{25} (2008), 298--304.

\bibitem{hoschek93}
Josef Hoschek and Dieter Lasser, \emph{{Fundamentals of Computer Aided
  Geometric Design}}, A. K. Peters, Ltd., USA, 1993.

\bibitem{jssz2013}
Gašper Jakli\v{c}, Maria~Lucia Sampoli, Alessandra Sestini, and Emil Žagar,
  \emph{{$C^1$ rational interpolation of spherical motions with rational
  rotation-minimizing directed frames}}, \CAGD \textbf{30} (2013), 159--173.

\bibitem{juettler99b}
Bert J\"uttler and Christoph M\"aurer, \emph{{Rational approximation of
  rotation minimizing frames usign Pythagorean hodograph cubics}}, \JGG
  \textbf{3} (1999), 141--159.

\bibitem{kps24}
Marjeta Knez, Francesca Pelosi, and Maria~Lucia Sampoli, \emph{{Construction of
  $G^2$ spatial interpolants with prescribed arc lengths}}, \JCAM \textbf{441}
  (2024), 115684.

\bibitem{ks21}
Marjeta Knez and Maria~Lucia Sampoli, \emph{{Geometric interpolation of ER
  frames with $G^2$ {Pythagorean}-hodograph curves of degree 7}}, \CAGD
  \textbf{88} (2021), 102001.

\bibitem{kv12}
Marjeta Krajnc and Vito Vitrih, \emph{{Motion design with Euler–Rodrigues
  frames of quintic {Pythagorean}-hodograph curves}}, \MCS \textbf{82} (2012),
  1696--1711.

\bibitem{wang97}
Wenping Wang and Barry Joe, \emph{Robust computation of the rotation minimizing
  frame for sweep surface modeling}, \CAD \textbf{29} (1997), 379--391.

\bibitem{wang08}
Wenping Wang, Bert Jüttler, Dayue Zheng, and Yang Liu, \emph{Computation of
  rotation minimizing frames}, \ACMTG \textbf{27} (2008), 1--18.

\bibitem{sj2007}
Zbyněk Šír and Bert Jüttler, \emph{{$C^2$ Hermite interpolation by
  {Pythagorean Hodograph} space curves}}, \MC \textbf{76} (2007), 1373--1391.

\end{thebibliography}

\end{document}